\documentclass{article}
\usepackage{PRIMEarxiv}
\usepackage{multirow}
\usepackage{graphicx} 
\usepackage{algorithm}
\usepackage{algorithmic}
\usepackage{amsmath}
\usepackage{amsfonts}
\usepackage{amssymb}
\usepackage{amsbsy}
\usepackage{bbm}
\usepackage{xcolor}
\usepackage{hyperref}
\newcommand{\E}{{\mathbb E}}
\newcommand{\calK}{\mathcal{K}}
\newcommand{\calR}{\mathcal{R}}

\newcommand{\lpar}{\left(}
\newcommand{\rpar}{\right)}

\usepackage[utf8]{inputenc} 
\usepackage[T1]{fontenc}    
\usepackage{hyperref}       
\usepackage{url}            
\usepackage{booktabs}       
\usepackage{amsfonts}       
\usepackage{nicefrac}       
\usepackage{microtype}      
\usepackage{lipsum}
\usepackage{fancyhdr}       
\usepackage{graphicx}       
\graphicspath{{media/}}     
\newtheorem{assumption}{Assumption}
\newtheorem{lemma}{Lemma}
\newtheorem{theorem}{Theorem}
\title{Fixed Confidence and Fixed Tolerance  Bi-level Optimization for Selecting the Best Optimized System
}

\author{
  Yuhao Wang \\
  School of Industrial and Systems Engineering\\
  Georgia Institute of Technology \\
  Atlanta\\
  \texttt{yuhaowang@gatech.edu} \\
  \And
 Seong-Hee Kim\\
  School of Industrial and Systems Engineering\\
  Georgia Institute of Technology \\
  Atlanta\\
  \texttt{skim@isye.gatech.edu}
   \And
  Enlu Zhou \\
  School of Industrial and Systems Engineering\\
  Georgia Institute of Technology \\
  Atlanta\\
  \texttt{enlu.zhou@isye.gatech.edu} \\
}

\begin{document}
\maketitle

\begin{abstract}
In this paper, we study a fixed-confidence, fixed-tolerance formulation of a class of stochastic bi-level optimization problems, where the upper-level problem selects from a finite set of systems based on a performance metric, and the lower-level problem optimizes continuous decision variables for each system. Notably, the objective functions for the upper and lower levels can differ. This class of problems has a wide range of applications, including model selection, ranking and selection under input uncertainty, and optimal design. To address this, we propose a multi-stage \textbf{Pruning-Optimization} framework that alternates between comparing the performance of different systems (\textbf{Pruning}) and optimizing systems (\textbf{Optimization}). 
This multi-stage framework is designed to enhance computational efficiency by pruning inferior systems with high tolerance early on, thereby avoiding unnecessary computational efforts. We demonstrate the effectiveness of the proposed algorithm through both theoretical analysis of statistical validity and sample complexity and numerical experiments.

\end{abstract}


\section{INTRODUCTION}
\label{sec:intro}



We study a problem of selecting the best among a finite set of systems, where the (expected) performance of each system is a  function of some decision variable. The best system is defined as the system that has the best expected performance under the ``optimized" decision variable. This optimized decision variable optimizes an objective that can be different from the system's expected performance. Mathematically, the problem can be formulated as a bilevel optimization of  the following form:
 
    \begin{align}
     \min_{k\in \mathcal{K}} \quad  &h(k,x^*_k) = \mathbb{E}_\zeta\left[H(k,x^*_k,\zeta)\right],  \label{eq: bilevel intro}
     \\
     s.t. ~~& x^*_k = \arg\min_{x_k\in\mathcal{X}_k} f(k,x_k) = \mathbb{E}_\xi \left[F(k,x_k,\xi)\right] \label{eq: bilevel intro2}
     \end{align} 
Here $\mathcal{K} =\{1,2,\ldots,K\}$ is the finite candidate set of systems, the upper-level objective $h(k,x_k^*)$ is the expected performance for system $k\in\mathcal{K}$, and $x_k^*$ is its optimized decision variable that optimizes the lower-level objective $f(k,x_k)$.  Given a fixed value of $x_k, k\in\mathcal{K}$, the upper-level problem \eqref{eq: bilevel intro} belongs to the regime of ordinal optimization with categorical decision variables; whereas the lower-level problem \eqref{eq: bilevel intro2} belongs to the regime of (continuous) stochastic optimization. This formulation covers a wide range of problems, such as those we list below.
\begin{enumerate}
 \item {\bf Model selection and parameter learning:} In machine learning problems, one often faces the problem of learning the model parameter using training data and selecting the best model (with the learned parameter) that minimizes the loss over the validation data. Suppose we have a training dataset $\{d_i=(x_i,h_i)\}_{i=1}^N$ and a validation dataset $\{v_i=(x_j,h_j)\}_{j=1}^M$. Given the $k^{th}$ model, its parameter is often learned by empirical risk minimization (ERM) given as follows
    \begin{equation}\label{ERM}
    \theta_k^* = \arg\min_{\theta}\sum_{i=1}^N l(k,\theta,d_i),
    \end{equation}
    where $l(k,\theta,\cdot)$ is a loss function over the training data, associated with the $k^{th}$ model and its parameter $\theta$. Then, the best model is selected by
    \begin{equation}\label{model_select}
    \min_{k=1,\ldots,K}\sum_{j=1}^N L(k,\theta_k^*,v_i),
    \end{equation}
    where $L(k,\theta_k^*,\cdot)$ is a loss function over the validation data, associated with the $k^{th}$ model and its learned parameter $\theta_k^*$. Putting in the form of \eqref{eq: bilevel intro} and \eqref{eq: bilevel intro2}, \eqref{ERM} is the lower-level optimization and \eqref{model_select} is the upper-level optimization.
    \item {\bf Ranking and Selection (R\&S) under input uncertainty:} In R\&S problems, one wants to select the best design that has the smallest mean performance among a finite set of designs $\{1,\ldots, K\}$, i.e., $\min_{k\in\{1,\ldots,K\}}\mathbb{E}_{\xi}[H(k,\xi)]$.  In practice, we only have a finite amount of data to estimate the input distribution of $\xi$. To hedge against the input uncertainty, \cite{zhou2015simulation} and \cite{ZhouWu:2018BRO} proposed the Bayesian Risk Optimization (BRO) approach. It quantifies the input uncertainty using the Bayesian posterior distribution on the unknown input parameter, assuming the input distribution takes a parametric form, and then imposes a risk measure with respect to the posterior on the unknown mean performance measure. That is, 
    \begin{equation}\label{BRO}
    \min_{k=1,\ldots,K} \rho_{\theta}\E_{\xi|\theta}[H(k,\xi)],
    \end{equation}
    where $\rho_{\theta}$ denotes the risk measure $\rho$ taken with respect to the posterior of the unknown input parameter $\theta$, and $\E_{\xi|\theta}$ denotes the expectation with respect to $\xi$ given a value of $\theta$. For coherent risk measures , it can be represented in the dual form (see, e.g., \cite{ang2018dual} for coherent risk measures and dual forms) that leads to an lower-level minimization. For example, if $\rho$ is set as the Conditional Value-at-Risk (CVaR) with risk level $\alpha$, then \eqref{BRO} can be re-written as
    \begin{equation}\label{eq: risk-averse R&S}
  \min_{k=1,\ldots,K}\min_{x\in\mathbb{R}} \left\{x + \frac{1}{1-\alpha}\E_{(\theta,\xi)}[(H(k,\xi)-x)_{+}]\right\},
    \end{equation}
    where $(\cdot)_+=\max(\cdot,0)$, $\E_{(\theta,\xi)}[\cdot] = \E_{\theta}\E_{\xi|\theta}[\cdot]$. Note \eqref{eq: risk-averse R&S} can be put in the form of \eqref{eq: bilevel intro} and \eqref{eq: bilevel intro2} with $H(k,x,\zeta)=F(k,x,\zeta) = \{x+\frac{1}{1-\alpha}(H(k,\xi)-x)_+\}$ and $\zeta$ replaced by $(\theta,\xi)$.
   \item {\bf Optimal design:} In problems of optimal design such as structural optimization and optimal shape design (e.g., \cite{herskovits2000contact}), in the upper level the decision maker needs to decide types of materials and shapes to minimize the total cost. For each structure, in the lower level the amount of materials needs to be decided by solving a potential energy minimization problem for achieving the best stability of the structure. Another example is optimal drug design (e.g., \cite{verweij2020randomized}), where a drug company needs to select among multiple possible drugs or treatments for deceases. The lower-level problem is to optimize the treatment effect for a specific drug/treatment by optimizing through variables such as the dosage of the drug. The upper-level problem is to select among possible drugs/treatments under its optimal dosage, possibly according to a different objective that balances the trade-off between cost and effectiveness.
\end{enumerate}

Motivated by these examples, we propose a \textbf{fixed confidence} and \textbf{fixed tolerance} (FCFT) formulation for finding a near  optimal system $\bar{k}$ of \eqref{eq: bilevel intro}. The FCFT formulation aims to select a system  $\Bar{k}$ whose optimal expected performance is no worse 
than that of the best system by a pre-specified \textbf{tolerance}, with a probability at least of a pre-specified \textbf{confidence} level. 
A valid FCFT procedure guarantees precision (confidence level and tolerance) for the final solution by either pre-specifying a computational budget or terminating based on a carefully designed stopping criterion. To our best knowledge, the FCFT formulation is new in the literature of bilevel optimization (see the survey paper \cite{sinha2017review}), as most of existing works on bilevel optimization seek to develop methods that converge fast but lack statistical guarantees.

On a related note, a similar criterion to FCFT has been recently proposed for R\&S problems (e.g. see \cite{eckman2021fixed}) to select a ``good system" with high probability. However, in R\&S the performance of each system is simply modelled as a random variable that does not depend on any decision variable and, the phrase ``good system" refers to systems whose expected performance is within the tolerance compared with the performance of the optimal system. 
Moreover, the concept of FCFT is also closely related to the probably approximately correct (PAC) learning, where one studies the sample complexity of a learning algorithm to produce a hypothesis that is approximately correct (i.e., within an small tolerance) with high probability (confidence level). However, the research on PAC learning focuses more on the rate of sample complexity rather than the exact number of samples with given tolerance and confidence level, which is the focus of our study. For implementation of an algorithm, only having the rate of sample complexity is not sufficient to compute an exact number of samples or design a stopping criteria that guarantee the required tolerance and confidence level. 

In this paper, we develop a novel multi-stage pruning-optimization framework that alternates between optimizing the decision variable for each remaining system and comparing the performance of different systems across various stages. 
In earlier stages, we solve both the upper-level and lower-level problems with larger tolerance value. By identifying and pruning systems that are sufficiently inferior to other systems, the computational effort is saved for later stages on systems that are more likely to be the optimal system. We summarize our main contributions as follows:
\begin{enumerate}
    \item We propose a new and practical problem formulation of FCFT bi-level optimization of selecting the best optimized system. To solve this problem, we propose a novel multi-stage \textbf{pruning}-\textbf{optimization} framework that alternates between solving the upper-level problem through a \textbf{pruning} step and solving the lower-level problem through an \textbf{optimization} step across stages.
    By utilizing a sequence of decreasing tolerance value, this framework saves computational effort by identifying and removing inferior systems early on.
    \item In the \textbf{pruning} step, inspired by a feasibility check procedure by \cite{zhou2022finding}, we develop a sequential algorithm that sequentially checks whether a system can be evidently identified as a sub-optimal system, under the current lower-level solutions. 
    \item In the \textbf{optimization} step, we first adopt the stochastic accelerated gradient descent algorithm and an existing concentration result (\cite{Lan:book2020}) to determine the total number of iterations that meets the FCFT requirement. While this is a valid approach, it is computationally inefficient from both of our theoretical and numerical analyses. To improve the computational efficiency, we establish the weak convergence of the function error of the upper-level objective evaluated under solutions solved in the lower-level problems using stochastic gradient descent (SGD), which can be of independent interest. The convergence results are different in cases with the same or different objectives for the upper and lower level. In both cases, we then determine the total number of SGD iterations by an asymptotic concentration bound based on the aforementioned weak convergence result. With this asymptotic approach, we show in both cases the computational complexity is significantly improved. 

    \item We demonstrate the effectiveness of the proposed algorithm both theoretically and numerically. Theoretically, we prove the statistical validity of the proposed algorithm, meaning the requirements on the confidence level and tolerance level are satisfied. Numerically, on an example of drug selection with optimal dosage, we show the multi-stage framework can significantly save computational effort compared with a benchmark approach.
\end{enumerate}


\subsection*{Literature Review}
\subsubsection*{Bi-level optimzation}

Bi-level optimization is a hierarchical optimization framework where one problem (the ``lower-level" problem) is nested within another (the ``upper-level" problem). The upper-level problem depends on the solution of the lower-level problem, making the process complex and computationally demanding. This structure is useful for real-world applications like machine learning, game theory, and economics, where decisions must be optimized in two interconnected stages. While bilevel optimization can model these hierarchical interactions effectively, solving such problems often requires specialized algorithms due to their non-convex nature and computational challenges. Many different approaches have been proposed for solving bilevel optimization. Classic approaches includes single-level reduction for linear bilevel optimization \cite{bard1982explicit,fortuny1981representation}, descent method \cite{kolstad1990derivative,savard1994steepest,vicente1994descent}, penalty function methods \cite{aiyoshi1981hierarchical,aiyoshi1984solution,ishizuka1992double} and trust-region methods \cite{liu1998trust,marcotte2001trust}. Later evolutionary algorithms are adopted to handle bi-level optimization \cite{mathieu1994genetic,yin2000genetic,hejazi2002linear,wang2005evolutionary,wang2006review}. All these aforementioned methods assumed the decision variables in both levels are continuous. There were also a stream of works studying discrete bilevel optimization, e.g., \cite{moore1990mixed,vicente1996discrete,bard1992algorithm}. We refer to \cite{sinha2017review} for a comprehensive overview on bilevel optimization. Methods for both continuous and discrete bilevel optimization utilize the inherent relations (such as distance, ordering) between different solutions, which does not exist for categorical decision variables. As a result, These methods cannot be efficiently applied to our setting. Furthermore, most of the aforementioned methods are designed for deterministic optimization, whereas our problem setting has stochastic optimization in both levels.

\subsubsection*{Ranking and selection (R\&S)}
Another closely related literature is R\&S, which aims to select the best among a finite set of designs. Classic R\&S (e.g., \cite{chen2000simulation,KN:2001}) can be viewed as a special case of our problem with the domain space of lower-level problem to be a singleton. 
R\&S has been studied in various different settings, which include but are not limited to contextual R\&S \cite{shen2021ranking,cakmak2023contextual,du2024contextual}, robust R\&S \cite{gao2017robust,fan2020distributionally},
constrained R\&S \cite{hunter2013optimal,pasupathy2014stochastically}, R\&S with input uncertainty \cite{corlu2015subset,wu2017ranking,xiao2018simulation,song2019input,xiao2020optimal,xu2020joint}, and data-driven R\&S \cite{wu2019fixed,wu2022data,wang2022fixed,kim2022optimizing,wang2023input,wang2024optimal}. Despite these different variants, they share in common that the performance of each design is simply modeled by a random variable that does not involve any extra decision variables. There is only one exception \cite{si2022selecting} that considers a similar setting as this paper, but with two major differences. First, they consider the same objective function in the upper and lower levels, whereas we consider a more general setting where the objectives can be different in the upper and lower levels. Second, they use a fixed budget formulation, meaning a total computing budget is pre-specified and the task is to allocate the budget efficiently. In contrast, our formulation takes a pre-specified confidence level and tolerance level as inputs, and our goal is to determine the computing budget to achieve the specified confidence and tolerance level. Because of the difference in formulations, the methodology we develop also diverges considerably from that of \cite{si2022selecting}.

\section{Problem Statement}
\label{sec:problem}
For ease of reading, we re-state Problem \eqref{eq: bilevel intro} and \eqref{eq: bilevel intro2} here.
\begin{equation} \label{eq: RSOS general}
    \begin{aligned}
     \min_{k\in \mathcal{K}} \quad  &h(k,x^*_k) = \mathbb{E}_{\zeta}\left[H(k,x^*_k,\zeta)\right], 
     \\
     s.t. ~~& x^*_k = \arg\min_{x_k\in\mathcal{X}_k} f(k,x_k) = \mathbb{E}_\xi \left[F(k,x_k,\xi)\right]
     \end{aligned} 
\end{equation}
where $\mathcal{K}=\{1,2,\ldots,K\}$ is a finite set of systems, $x_k \in \mathcal{X}_k \subseteq \mathbb{R}^{d_k}$ is a continuous decision variable that affects the lower-level (expected) performance $f(k,x_k)$ and upper-level (expected) performance $h(k,x_k)$ for each $k \in \calK$, and $\mathcal{X}_k$ is the feasible set of $x_k$. Let $k^* = \arg\min_{k\in \mathcal{K}} h(k,x^*_k)$ denote the optimal system.

In the lower level, for a fixed system $k$ one  can utilize an optimization algorithm to solve $f(k,\cdot)$ for a near-optimal solution that approximates $x_k^*$. In the upper level, for a fixed system $k$ and a fixed decision variable $x_k$, one has access to a sequence of independent and identically distributed (i.i.d.) random function evaluations $H(k,x_k,\zeta_1), H(k,x_k,\zeta_2),\ldots$ with mean $h(k,x_k)$. These random samples can then be used to obtain an estimate of $h(k,x_k)$. Due to the randomness in both stochastic gradients and function evaluations, it is impossible to find the true optimal system $k^*$ 
with probability $1$ under finite samples. As a consequence, we want to achieve the following goal:
with probability at least $1-\alpha$, find $\bar{k}$ such that 
$$ h(\Bar{k}, x_{\Bar{k}}^*) \le  h(k^*,x^*_{k^*}) + \epsilon,$$
where $\epsilon>0$ represents the \textbf{error tolerance} of the decision maker. We say $\Bar{k}$ is an $\epsilon$-optimal system with \textbf{confidence level} $1-\alpha$.

\section{Multi-stage Framework}
\label{sec: framework}

\subsection{A Direct Two-Step Optimization-Evaluation Approach} \label{sec: direct approach}
To address the problem outlined in \eqref{eq: RSOS general}, a straightforward strategy is to initially identify a near-optimal decision variable \({x}_k\) for each system \(k \in \calK\), using an optimization algorithm such as stochastic gradient descent (SGD). Subsequently, for each system \(k \in \calK\), an estimator \(\widehat{h}_k\) of its expected upper-level performance \(h(k,{x}_k)\) under the decision variable \({x}_k\) is obtained through repeated function evaluations. 
Suppose that, for lower-level problem, we achieve error tolerance $\epsilon_1$ with confidence level $1-\alpha_1$, i.e.,
with probability at least $1-\alpha_1$, ${x}_k$ satisfies $h(k,{x}_k) - h(k,x_k^*) \le \epsilon_1, \forall k\in \calK$. Furthermore, suppose that, for upper-level comparison, the pair-wise performance difference evaluation achieves $\epsilon_1'$ error tolerance with confidence level $1-\alpha_1'$, i.e., with probability at least $1-\alpha'_1$, $|(\widehat{h}_k-\widehat{h}_i) - \lpar h(k,{x}_k)-h(i,{x}_i)\rpar |\le \epsilon'_1$ for each $i,k \in\mathcal{K}$. 
Then, if we pick $\Bar{k}$ such that
$\Bar{k} = \arg\min_k \widehat{h}_k $, we have with probability at least $1-\alpha_1-\alpha'_1$ (by union bound), 
\begin{align*}
    h(\Bar{k},x_{\Bar{k}}^*) 
    &\le h(\Bar{k},{x}_{\Bar{k}}) \\
    &\le h(k^*,{x}_{k^*}) + \widehat{h}_{\Bar{k}}-\widehat{h}_{k^*} + \epsilon'_1  \\
    &\le h(k^*,x_{k^*}^*)+ \widehat{h}_{\Bar{k}}-\widehat{h}_{k^*}+\epsilon_1 + \epsilon'_1 \\
    &\le h(k^*,x_{k^*}^*)+ \epsilon_1 + \epsilon'_1 . 
\end{align*}
That is, $\Bar{k}$ is an $(\epsilon_1+\epsilon'_1)$-optimal  system with confidence level $1-\alpha_1-\alpha'_1$. By choosing $\alpha_1+\alpha'_1 = \alpha$ and $\epsilon_1+\epsilon'_1 = \epsilon$, we obtain that $\Bar{k}$ is an $\epsilon$-optimal system. 

While this described methodology achieves our objective, its practical application may suffer from computational inefficiency. Specifically, solving all systems to optimality or estimating their expected performance within narrow error tolerance (here, $\epsilon_1$) requires substantial computational efforts. When there is a large performance difference between systems, exhaustive optimization across all systems or precise estimation of their expected performances becomes unnecessary. Consider the following illustrative example with a set of two system \(\calK=\{1,2\}\) with expected performances defined by \(h(1,x_1) = f(1,x_1) = x_1^2 + 10\) and \(h(2,x_2) = f(2,x_2) = x_2^2\). Suppose the objective is to identify a near-optimal system $\Bar{k}$ with tolerance $\epsilon=1$, i.e., $\Bar{k} = 2$, with a confidence level of at least \(1-\alpha\). Employing a direct approach necessitates determining $\epsilon_1$ and $\epsilon'_1$ such that $\epsilon_1+\epsilon'_1 = \epsilon = 1$ and solving both $f(1,x_1)$ and $f(2,x_2)$ to  $\epsilon_1\le 1$-optimal. 
However, by selecting \(\epsilon_1 = 5\) and \(\epsilon'_1 = 4\), it can be demonstrated that, with a probability of at least \(1-\alpha\):\\
\centerline{
\(
\widehat{f}_2 - \widehat{f}_1 \le f(2,\bar{x}_2) - f(1,\bar{x}_1) + \epsilon'_1 \le f(2,x_2^*) - f(1,x_1^*) + \epsilon_1 + \epsilon'_1 = 0 - 10 + 5 + 4 = -1 < 0.
\)
}
\\
This identifies system 2 as the optimal system despite the larger tolerances \(\epsilon_1\) and \(\epsilon'_1\). Once we identify system 2 as the optimal system with $\epsilon_1= 5$ and $\epsilon_1' = 4$, we can spend the remaining computational resources to optimally solve only system 2 to the desired tolerance. 
This approach with larger $\epsilon,\epsilon'$ enables a more efficient allocation of computational resources by focusing efforts on systems with a higher likelihood of being optimal. Motivated by this, in the next section we introduce a novel multi-stage pruning-optimization framework.

\subsection{Multi-Stage Pruning-Optimization Framework} 

Addressing the issue of  sampling inefficiency inherent in the aforementioned direct approach, we propose a multi-stage pruning-optimization framework that utilizes a sequence of progressively narrowing error tolerances. Specifically, the framework operates on two sequences of tolerances:  $\epsilon_1 > \epsilon_2 >\ldots>\epsilon_T$ and $\epsilon'_1>\epsilon'_2>\ldots > \epsilon'_T$, and two sequences of significance levels (significant level is defined as one minus the confidence level): $\{\alpha_t\}_{t=1}^T $ and $\{\alpha'_t\}_{t=1}^T $.  Here $\epsilon_t$ and $\alpha_t$ are the tolerance and significance level for solving $h(k,x_k^t)$ with approximate solution ${x}_k^t$ at stage $t$, respectively; whereas $\epsilon'_t$ and $\alpha'_t$ are the tolerance and significance level for estimating the performance difference between system $k$ and $i$ at stage $t$, respectively. We require $\epsilon_T$ and $\epsilon'_T$ satisfy $\epsilon'_T + \epsilon_T = \epsilon$ and $\{\alpha_t,\alpha'_t\}_{t=1}^T$ satisfy $\sum_{t=1}^T \alpha_t + \sum_{t=1}^T \alpha'_t = \alpha$. 
The methodology begins with a comparatively larger tolerance value, which is then gradually narrowed down to the target tolerance. At any given stage \(t \leq T\), the process involves optimizing each of the "remaining systems" (i.e., those not yet identified as sub-optimal) to find a solution that is \(\epsilon_t \)-optimal.
By discarding systems identified to be inferior relative to others early on, computation resources can be saved for later stages. This multi-stage pruning-optimization process efficiently identifies the optimal system by progressively concentrating efforts on the most promising candidates, thereby enhancing overall efficiency. 

The number of tolerance levels, \(T\), plays the role of balancing computational efficiency and accuracy of the optimization and simulation process. With a larger $T$, larger tolerances are adopted in early stages, and computational effort can be potentially saved by removing inferior systems with difference from any other system within larger tolerances. However, a larger $T$ also requires division of the total significance level $\alpha$ into more stage-wise significance levels $\alpha_t,\alpha'_t$, which can potentially increase the computational cost for large $T$. On a related note, when $T=1$, the procedure reduces to the direct approach.  In practice, one possible choice for the number of tolerance and tolerance sequence is \(T = \lceil \log_2 K \rceil\), where \(K\) represents the total number of systems, and $\epsilon'_t + \epsilon_t = 2^{T-t-1}\epsilon,~t=1,\ldots,T$.  In this way, 
when there is a large number of systems whose best performance are approximately uniformly distributed, every time by halving the tolerance value, one expects half of the remaining sub-optimal systems can be identified and removed.


To be specific, let $\calR_t$ be the remaining set at the beginning of stage $t$, with $\calR_1 = \calK$ representing the initial set of all systems.   
 At the start of each stage \(t\), with initial solutions \({x}_k^{t-1}\) for each system \(k\) within \(\calR_t\), a predetermined number  of SGD steps are executed for system \(k\). This process yields updated solutions \({x}_k^t\) for all \(k \in \calR_t\), ensuring that, with probability \(1 - \alpha_t\),

 \centerline{$ h(k,{x}_k^t) \le h(k,x_k^*) + \epsilon_t, \quad \forall k \in \calR_t.$}
\noindent This procedure is referred to as the \textbf{Optimization} step at stage \(t\), which solves the lower-level problem to an intermediate tolerance $\epsilon_t$.

Following the {\bf Optimization} step at stage \(t\), in the upper level, we next construct an estimator \(\widehat{h}_k^t-\widehat{h}_i^t\) that approximates the difference \(h(k,{x}_k^t) - h(i,{x}_i^t)\) for all pairs of systems \(i \neq k\) within the remaining set \(\calR_t\). This estimator is constructed by executing a number of simulation runs, which  guarantee the estimator \(\widehat{h}_k^t-\widehat{h}_i^t\) achieves an error tolerance of \(\epsilon'_t\) with a specified probability \(1-\alpha'_t\).



With a combined probability of \(1 - \alpha_{t} - \alpha'_{t}\) by the union bound, it follows that for every pair of distinct systems \(k\) and \(i\) within the remaining set \(\calR_t\),
$ h(k,x_k^*) - h(i,x_i^*) \ge h(k,{x}_k^t) - h(i,{x}_i^t) - \epsilon_t \ge \widehat{h}_k^t - \widehat{h}_i^t - \epsilon_t - \epsilon'_t$, 
for all \(k \neq i \in \calR_t\). Therefore, if for any system \(i \in \calR_t\), the condition \(\widehat{h}_k^t - \widehat{h}_i^t \ge \epsilon_t + \epsilon'_t\) is met, system \(k\) can be eliminated from \(\calR_t\) as it is no better than system \(i\). In the case where no such system \(i\) meets this criterion with respect to system \(k\), then system \(k\) is retained within \(\calR_t\). This process is referred to as the \textbf{Pruning} step in stage \(t\).



The multi-stage \textbf{Pruning-Optimization} Framework is outlined in Algorithm \ref{alg: RSCOP}. 
Details regarding the \textbf{Pruning} and \textbf{Optimization} steps are elaborated in Sections  \ref{sec: pruning} and \ref{sec: optimization}, respectively. 

\begin{algorithm}[H]
   \caption{Multi-stage Pruning-Optimization Framework}
   \label{alg: RSCOP}
\begin{algorithmic}[1]
   \STATE {\bfseries Input:} System set $\calK$, sequences of tolerance ($\epsilon_1,\ldots,\epsilon_T $) and ($\epsilon'_1,\ldots,\epsilon'_T$), initial point for decision variable ${x}_k^0$ for $k\in \calK$, input parameter $r_0$ for \textbf{Pruning}.
   
   \STATE { \bfseries Initialize:}  $\calR_0 \leftarrow \calK$. $t\leftarrow1$.
   \WHILE{$t < N$ and $|\calR_{t}| > 1$}
   \STATE For each $k\in \calR_t$, update the decision variable to ${x}_k^t$ by running \textbf{Optimization}. 
   \STATE  Update the remaining set $\calR_{t+1}$ by running \textbf{Pruning}. 
   \STATE $t\leftarrow t+1$.
   \ENDWHILE
   \STATE {\bfseries Output:} Any $k\in \calR_t$ as an $\epsilon$-optimal system.
\end{algorithmic}
\end{algorithm}

\section{Fully Sequential Pruning Step} \label{sec: pruning}

In this section, we elaborate on the {\bf Pruning} step introduced in Section \ref{sec: framework}. Recall in Section \ref{sec: framework} we want to construct estimator $\widehat{h}_k^t, \forall k\in \calR_t$ such that 
$|\widehat{h}_k^t - \widehat{h}_i^t - ((h(k,{x}_k^t) - h(i,{x}_i^t)) | \le \epsilon'_t, \forall k\neq i \in \calR_t$.
When i.i.d. samples of function evaluation are available, such estimator can be obtained by sample average estimator with sample size determined by leveraging some concentration inequality of $H(k,{x}_k^t,\xi), k\in\calR_t$. For notation simplicity, in the remaining of the section, let $h^t_k = h(k,{x}^t_k), \forall k \in \calK$ be the expected outcome for system $k$ at decision variable ${x}_k^t$, and let $H^t_k = H(k,{x}^t_k,\xi)$ be corresponding stochastic function evaluation.

Nonetheless, such a static approach that pre-determines the number of simulations at the beginning can be inefficient in terms of the total number of function evaluations, especially when compared to a fully sequential approach that  determines whether more simulation is needed sequentially. 
Specifically, the reason of constructing the estimator $\widehat{h}_k^t-\widehat{h}_i^t$ is to infer the value of $h_k^t-h_i^t$. We are interested  in whether the following two inequalities are satisfied: for any two systems $i\neq k$ and $i,k \in \calR_t$,
\begin{equation}\label{eq: eliminate i}
    h_i^t - h_k^t \le \epsilon'_t 
\end{equation}
\begin{equation} \label{eq: eliminate k} 
        h_i^t - h_k^t \ge -\epsilon'_t 
\end{equation}
To see why \eqref{eq: eliminate i} and \eqref{eq: eliminate k} are of interest, if there exists $i\neq k$  such that \eqref{eq: eliminate i} does not hold, then we know with probability $1 - \alpha_{t}$,
\begin{equation} \label{eq: eliminate i reason}
\begin{aligned}
    h(i,x^*_i) &\ge h_i^t - \epsilon_t \\
    &>  h_k^t + \epsilon'_t - \epsilon_t 
    \\&\ge h(k,x^*_k) + \epsilon'_t - \epsilon_t.
\end{aligned}
\end{equation}
The first inequality holds since $x_k^t$ is an $\epsilon_t$-optimal solution solved in the optimization step.
  Suppose $\epsilon'_t > \epsilon_t$, which can be easily satisfied as the tolerance sequences are chosen by the decision maker. Then we obtain $h(i,x^*_i) > h(k,x^*_k)$,  which implies system $k$ is superior to system $i$ with probability $1-\alpha_t$. Hence, there is no need to further optimize or evaluate system $i$, and we remove system i from the remaining set, i.e., let $\calR_t = \calR_t \backslash \{i\}$. Otherwise, if \eqref{eq: eliminate i} holds for a fixed $i$ and all $k\in\calR_t$, we obtain
 $$ h_i^t-h(k,x_k^*) \le h_i^t - h_k^t + \epsilon_t \le \epsilon'_t+ \epsilon_t.$$
 Setting $k=k^*$, we have that $i$ is an $(\epsilon_t+\epsilon'_t)$-optimal system.
Similarly, if there exists $i\neq k$ such that \eqref{eq: eliminate k} does not hold, then $k$ is identified as an inferior system and we let $\calR_t = \calR_t \backslash \{k\}$. Otherwise, if for a fixed $k$ and all $i\in\calR_t$ such that \eqref{eq: eliminate k} holds, we have $k$ is an  $(\epsilon_t+\epsilon'_t)$-optimal system. In summary, if a pair of systems, $i$ and $k$, satisfy both \eqref{eq: eliminate i} and \eqref{eq: eliminate k}, then both systems remain ${\cal R}_t$. However, if either \eqref{eq: eliminate i} or \eqref{eq: eliminate k} is infeasible, one system is pruned.

To efficiently determine the feasibility with respect to \eqref{eq: eliminate i} and \eqref{eq: eliminate k} in a sequential manner, we draw upon the methodology outlined in \cite{zhou2022finding},  which provides a fully sequential feasibility check procedure to determine the feasibility of constraints.  We construct a confidence interval of the target value  $h_i^t - h_k^t$, whose width depends on the confidence level and a parameter $\tau_t$ chosen by the decision maker. The confidence interval shrinks as more samples are available and possesses the following property: with high probability, the target value can exit this interval by at most $\tau_t$ throughout the entire sampling process, where $\tau_t$ will be specified later.  The reason of introducing $\tau_t$ is to ensure the algorithm will terminate within finite steps.  Moreover, a larger $\tau_t$ results in a narrower confidence interval, which can help reduce the computational effort. However, it also leads to under coverage, i.e., the target value is not contained in the confidence interval. During the procedure, we will monitor whether the upper or lower bound of the constructed confidence interval hits a threshold value $q_t$. 
When the upper bound of the confidence interval reaches \(q_t\), the procedure declares \eqref{eq: eliminate i} as feasible and with high probability, \(h_i^t - h_k^t \le q_t + \tau_t\). 
To ensure feasibility, we want \eqref{eq: eliminate i} to be declared feasible only if 
\[
h_i^t - h_k^t \leq \epsilon'_t.
\]
This condition can be satisfied by setting \(q_t + \tau_t = \epsilon'_t\).
Conversely, if the lower bound of the confidence interval reaches \(q_t\), \eqref{eq: eliminate i} is declared infeasible. In this case, with high probability, it holds that 
\[
h_i^t - h_k^t \geq q_t - \tau_t.
\]
By following a similar reasoning as in \eqref{eq: eliminate i reason}, we can deduce that 
\[
h(i, x^*_i) \geq h(k, x^*_k) + q_t - \tau_t - \epsilon_t.
\]
When \eqref{eq: eliminate i} is declared infeasible, we identify \(i\) as an inferior system and remove it. To ensure that \(i\) is indeed an inferior system with high probability, we set \(q_t - \tau_t = \epsilon_t\), which leads to 
\[
h(i, x^*_i) \geq h(k, x^*_k).
\]
By combining \(q_t + \tau_t = \epsilon'_t\) and \(q_t - \tau_t = \epsilon_t\), we derive 
\[
q_t = \frac{\epsilon'_t + \epsilon_t}{2}, \quad \tau_t = \frac{\epsilon'_t - \epsilon_t}{2}.
\]
The feasibility of \eqref{eq: eliminate k} can be determined in a similar manner.

Based on the outcome of this feasibility check, we can then proceed with the {\bf Pruning} step as discussed previously, effectively eliminating inferior systems from consideration.

 Let $H_{k,1}^t,H_{k,2}^t,\ldots$ denote a sequence  of samples of $H_k^t$ at current solution ${{x}^t_k}$ for all $k\in\mathcal{K}$. We make the following assumption for the validity of the \textbf{Pruning} step:\\

\begin{assumption} \ \label{assump: prunning}
\begin{enumerate}
    \item The random simulation outputs follow a multivariate normal distribution:
    \begin{equation}
        \left[
        \begin{aligned}
            H_1^t\\
            H_2^t\\
            \vdots\\
            H_K^t
        \end{aligned}
        \right]
        \sim 
        \mathcal{N}\left(
        \left[
        \begin{aligned}
            h_1^t\\
            h_2^t\\
            \vdots\\
            h_K^t\\
        \end{aligned}
        \right],
        \Sigma(\mathbf{\Bar{x}}^t)
        \right),
    \end{equation}
    where $\mathbf{{x}}^t = ({x}_1^t,\ldots,{x}_K^t)$ and $\Sigma(\mathbf{{x}}^t)$ is a $K\times K$ covariance matrix of $\mathbf{{x}}^t$.
    \item Conditioned on $\calR_t$ and $\mathbf{\Bar{x}^t}$, $\{H_{i,\ell}^t-H_{k,\ell}^t\}_{\ell=1}^\infty$ are independent and identically distributed (i.i.d.) for $i < k\in\calR_t$.
    \item The tolerance sequence $\{\epsilon_t\}_{t=1}^T$ and $\{\epsilon'_t\}_{t=1}^T$ satisfy $\epsilon'_t > \epsilon_t, t=1,\ldots,T$.
\end{enumerate}
\end{assumption}
Assumption \ref{assump: prunning}.1 assumes Gaussian noises of function evaluation, which is a common assumption and can be asymptotically satisfied with batched data. Here we assume the more general joint normal distribution, which can potentially help reduce variance of $H_{i,\ell}^t-H_{k,\ell}^t$, the pair-wise difference of function evaluations, when correlated samples are used. For example, when the function  is evaluated through simulation, the technique of common random number can be employed in generating $H_{i,\ell}^t-H_{k,\ell}^t$ to reduce its variance and saves computational effort.  Assumption \ref{assump: prunning}.2 assumes i.i.d.\@ samples of function evaluations, and Assumption \ref{assump: prunning}.3 can be easily satisfied as both $\{\epsilon_t\}_{t=1}^T$ and $\{\epsilon'_t\}_{t=1}^T$ are chosen by the decision maker. 

With Assumption \ref{assump: prunning}, the confidence interval of $h(i,{x}_i^t) - h(k,{x}_k^t)$ with $r$ samples is specified as $\Bar{H}_i - \Bar{H}_k \pm R(r;\tau_t,\eta,S_{ik}^2)/r$, where (i) $\Bar{H}_i$ and $\Bar{H}_k$, estimators of $h_i^t$ and $h_k^t$, are sample average of $r$ i.i.d.\@ simulation outputs from system $i$ and $k$, respectively; and (ii) $R(r;\tau_t,\eta,S_{ik}^2)$ is a non-negative function that decreases to 0 as $r$ increases, which is computed as 
$$
R\left(r; v, \eta, S_{ik}^2\right) := \max \left\{0, \frac{\left(r_0-1\right) \eta S_{ik}^2}{\tau_t}-\frac{\tau_t}{2} r\right\}.
$$ 
Here, $\eta = \frac{1}{2}\left(\left(\frac{2\alpha'_t}{|\calR_t|(|\calR_t|-1)}\right)^{-\frac{2}{r_0-1}}-1\right)$ is a function of the significance level $\alpha'_t := \frac{\alpha}{2N}$ for the feasibility check, and $z = S_{ik}^2$ is a one-time variance estimator for $H_i^t-H_k^t$ with $r_0$ samples.
We refer the reader to \cite{KN:Select2006} for more details on the construction of the confidence interval. we present the full description of the \textbf{Pruning} step
in Algorithm \ref{alg: pruning efficient}.

 \textbf{Remark:} Compared with the feasibility check in \cite{zhou2022finding}, Algorithm \ref{alg: pruning efficient} is more sample-efficient thanks to the special structure of the considered problem here. In \cite{zhou2022finding}, they need to generate samples of function evaluations individually for each tested constraints (but with multiple right-hand-side threshold values). In our context, \eqref{eq: eliminate i} and \eqref{eq: eliminate k} are treated as one constraint with two different threshold values. Since there are in total ${|\calR_t|(|\calR_t|-1)}/{2}$ different constraints and evaluating a single constraint one time requires two samples of function evaluation for two different systems, it requires a total number of $|\calR_t|(|\calR_t|-1)$ function evaluations to evaluate all constraints once. Notably, in our specific context, the function $h(k,x_k^t)$ appears in $|\calR_t|-1$ different constraints, indicating a more efficient way of generating samples by using the same sample of function evaluation in different constraints. This leads to a total number of only $|\calR_t|$ function evaluations for evaluating all constraints once.

\begin{algorithm}[h]
   \caption{ Fully Sequential Pruning}
   \label{alg: pruning efficient}
\begin{algorithmic}[1]
   \STATE {\bfseries Input:} Remaining set $\calR_t$; current solution ${x}_k^t$ for $k\in \calR_t$; tolerance $\epsilon_t,\epsilon'_t$; initial number of simulations $r_0$.
   \STATE { \bfseries Initialize:} 
    Set $q_t \leftarrow \frac{\epsilon_t + \epsilon'_t}{2}$ and $ \tau_t \leftarrow \frac{\epsilon'_t - \epsilon_t}{2}$; set $\alpha'_t \leftarrow \frac{\alpha}{2N}$;  compute $\eta \leftarrow \frac{1}{2}\left(\left(\frac{2\alpha'_t}{|\calR_t|(|\calR_t|-1)}\right)^{-\frac{2}{n_0-1}}-1\right)$; set $\operatorname{ON} \leftarrow \calR_t$; and set $\operatorname{STOP}_{ik,1} \leftarrow \operatorname{False}$ and $ \operatorname{STOP}_{ik,2} \leftarrow \operatorname{False}$, $ \forall i<k \in \calR_t$.
    \STATE  For each $k \in \calR_t$, generate $r_0$ i.i.d.\@ samples  $(H_{k,\ell})_{\ell=1}^{r_0}$ of $H_k^t$ and compute $\Bar{H}_{k} \leftarrow \frac{1}{r_0} \sum_{\ell=1}^{r_0} H_{k,\ell}$;
    compute $S^2_{ik} \leftarrow \frac{1}{r_0-1} \sum_{\ell=1}^{r_0} \left( H_{i,\ell} - H_{k,\ell} - \Bar{H}_{i}-\Bar{H}_{k} \right)^2$, $\forall i < k\in\calR_t$;  and set $r \leftarrow r_0$. \\[6pt]
   \WHILE{$|\operatorname{ON}| \ge 2$}
   \FOR{$i< k \in \operatorname{ON}$}
   \STATE $Z_{ik} \leftarrow R(r;\tau_t,\eta,S_{ik}^2)$
   \IF{$\operatorname{STOP}_{ik,1} = \operatorname{False}$ and $\Bar{H}_{i}-\Bar{H}_{k} - \frac{Z_{ik}}{r} \ge q_t$}
   \STATE $\operatorname{STOP}_{ik,1} = \operatorname{True}$, $\calR_t = \calR_t \backslash \{i\}$, $\operatorname{ON} = \operatorname{ON} \backslash \{i\}$.
   \ELSIF{$\operatorname{STOP}_{ik,1} = \operatorname{False}$ and $\Bar{H}_{i}-\Bar{H}_{k} + \frac{Z_{ik}}{r} \le q_t$}
   \STATE $\operatorname{STOP}_{ik,1} = \operatorname{True}$. 
   \ENDIF
   \IF{$\operatorname{STOP}_{ik,2} = \operatorname{False}$ and $\Bar{H}_{i}-\Bar{H}_{k} + \frac{Z_{ik}}{r} \le -q_t$}
   \STATE $\operatorname{STOP}_{ik,2} = \operatorname{True}$, $\calR_t = \calR_t \backslash \{k\}$, $\operatorname{ON} = \operatorname{ON} \backslash \{k\}$.
   \ELSIF{$\operatorname{STOP}_{ik,2} = \operatorname{False}$ and $\Bar{H}_{i}-\Bar{H}_{k} - \frac{Z_{ik}}{r} \ge -q_t$}
   \STATE $\operatorname{STOP}_{ik,2} = \operatorname{True}$.
   \ENDIF
   \ENDFOR
   \FOR{ $k \in \operatorname{ON}$}
   \IF{$ \operatorname{STOP}_{ik,1} = \operatorname{STOP}_{ik,2} = \operatorname{True}, \forall i<k \in \operatorname{ON}$ and $ \operatorname{STOP}_{ki,1} = \operatorname{STOP}_{ki,2} = \operatorname{True}, \forall i>k \in \operatorname{ON} $ }
   \STATE $\operatorname{ON} = \operatorname{ON} \backslash \{k\}$.
   \ENDIF
   \ENDFOR
   \STATE $r\leftarrow r+1$. Generate $1$ more sample $H_{k,r}$ of $H_{k}$ and update $\Bar{H}_{k}, \forall k \in \operatorname{ON}$.
   \ENDWHILE
   \STATE{\bfseries Output:} Set $\calR_t$.
\end{algorithmic}
\end{algorithm}

\section{Optimization Step} \label{sec: optimization}
Recall the goal of the \textbf{Optimization} step in Algorithm \ref{alg: RSCOP} is to solve a near-optimal solution ${x}_k^t$ such that $h(k,x_k^t)$ is within the tolerance $\epsilon_t$ compared with $h(k,x_k^*)$.  Generally, the difference between $h(k,{x}_k^t)$ and $h(k,x_k^*)$ can be bounded by difference between the near-optimal solution $x_k^t$ and the optimal solution $x_k^*$ with regularity condition such as Lipschitz continuity on $h(k,\cdot)$.
However , the near optimal solution ${x}_k^t$ is solved by optimizing the lower-level objective $f(k,\cdot)$ instead of the upper-level objective $h(k,\cdot)$. Bounding the error in the objective that we intend to optimize, which is $f(k,\cdot)$, usually requires convexity of $f(k,\cdot)$, whereas bounding the error in solution requires stronger assumption such as strong convexity. 
The following Assumption \ref{assump: convex} summarizes the regularity conditions on $f(k,\cdot)$ and $h(k,\cdot)$ for the purpose of solving the lower level problem to get an near-optimal solution $x_k^t$. 
\begin{assumption}\label{assump: Lipschitz-convex}
For each $k \in \mathcal{K}$,    
    \begin{enumerate}
        \item $f(k,x_k)$ is $\mu_k$-strongly convex in $x_k \in \mathcal{X}_k$. i.e., $f(k,y)-f(k,x)-\left\langle \nabla f(k,x), y-x\right\rangle \ge \frac{\mu_k}{2} \|x-y\|^2_2, \forall x,y \in \mathcal{X}_k$. 
        \item $h(k,x_k)$ is $L_k$-Lipschitz continuous in $x_k \in \mathcal{X}_k$.
    \end{enumerate}
\end{assumption}

With Assumption \ref{assump: Lipschitz-convex}, we can bound the error of $h(k,x_k)$ as 
$$|h(k,{x}_k) - h(k,x_k^*)| \le L_k \|{x}_k-x_k^*\|_2 \le L_k \sqrt{\frac{2(f(k,{x}_k) - f(k,x_k^*))}{\mu_k}}.$$
Consequently, if we can obtain an $x^t_k$ such that 
\begin{equation} \label{eq: near optimal f}
    f(k,x_k^t) - f(k,x_k^*) \le \frac{\epsilon_t^2 \mu_k}{2L_k^2},
\end{equation} $|h(k,x_k^t)-h(k,x_k^*)| \le \epsilon_t $, i.e., we reach the desired stage-wise tolerance.

To control the optimization error in \eqref{eq: near optimal f}, different algorithms and the associated finite-sample error bounds can be employed to determine the total number of iterations needed to achieve the desired tolerance and confidence level. Here in Section \ref{sec: optimization non-asymptotic} we take the stochastic accelerated descent (SAGD) algorithm from \cite{Lan:book2020} as an example and demonstrate how to design the \textbf{Optimization} step utilizing SAGD. However, such non-asymptotic approach is usually sample inefficient in practice, i.e., unnecessarily large number of iterations will be conducted, leading to much smaller optimization error than the desired tolerance. In practice, asymptotically valid concentration bound such as derived from central limit theorem is usually preferred to improve the sample efficiency. For this purpose, we develop an asymptotic approach in Section \ref{sec: optimization asymptotic} by studying the asymptotic convergence of the stochastic gradient descent (SGD) algorithm and show that the sample complexity is significantly improved using the asymptotic approach.

\subsection{Non-asymptotic Approach with Stochastic Accelerated Gradient Descent (SAGD)} \label{sec: optimization non-asymptotic}
In this section, we employ the stochastic accelerated gradient descent (SAGD) from \cite{Lan:book2020} as the optimization tool, which exploits the strong convexity of $f(k,\cdot)$. Since SAGD is not the focus of this paper, we defer the details of Assumptions required by SAGD and update rule of SAGD to the electronic companion Section \ref{ECsec: SAGD}. Furthermore, since we focus on optimizing a fixed system $k$ in the remaining of the section, we drop $k$ from the notation.

Let $x_N$ denote the solution given by SAGD after $N$ iterations. We have the following result from \cite{Lan:book2020} on the optimality gap.


\begin{lemma} \textbf{(Proposition 4.6 in \cite{Lan:book2020})} \label{lem: SAGD concentration}
Suppose Assumption \ref{assump: convex} holds. Let $({x}_{\ell})_{\ell\ge1}$ be the sequence given by SAGD.  Then, $\forall \lambda\ge0$, 
$$
\begin{aligned}
\mathbb{P}\left(f\left({x}_{{ N}}\right)-f(x^*) \geq   \frac{2 \lambda \sigma_{G} D }{\sqrt{{  3N}}} + \frac{4M^2 + 4(1+\lambda)\sigma_{G}^2}{\mu(N+1)}  + \frac{4\nu\|x_{0} - x^*\|_2^2}{N(N+1)t} \right)
 \le \mathrm{e}^{-\lambda} + \mathrm{e}^{-\frac{\lambda^2}{3}} .
\end{aligned}
$$
Here $x_0$ is the initial solution and $\sigma_G,D,M$ and $v$ are constant parameters specified in Section \ref{ECsec: SAGD}.
\end{lemma}

 The number  of SAGD iterations for a fixed system to guarantee $\epsilon_t$ tolerance (of $h$) with confidence level  $1-\alpha_t/K$ can be obtained by first calculating
\begin{equation} \label{eq: calculate lambda}
    \lambda_t = \min\left\{\lambda':  \mathrm{e}^{-\lambda'} + \mathrm{e}^{-\frac{(\lambda')^2}{3}} \le \frac{\alpha_t}{K} \right\},
\end{equation}
 and then setting
 \begin{equation} \label{eq: calculate N_k}
 \begin{aligned}
      N^t = \min&\left\{ N: \frac{2 \lambda_t \sigma_{G} D }{\sqrt{3N}} + \frac{4M^2 + 4(1+\lambda_t)\sigma_{G}^2}{\mu(N+1)} \right.\\
      &\left.~~~~~~~~+ \frac{4\nu D^2}{N(N+1)} \le \frac{\mu \epsilon^2_t }{2 L^2} \right\}
 \end{aligned}
 \end{equation}
Both \eqref{eq: calculate lambda} and \eqref{eq: calculate N_k} can be easily solved, for example, by binary line search.
The significance level \(\alpha_t\) is set to \(\frac{\alpha}{2T}\). 
Then, starting with the current solution ${x}_{N^{t-1}}$, one can conduct additional $N^t - N^{t-1}$ steps of SAGD to obtain ${x}^t = {x}_{N^{t}}$. Furthermore, \eqref{eq: calculate N_k} yields  a sample complexity of order, $N^t = O\left( \frac{\log^2\left(\frac{T}{\alpha}\right)}{\epsilon_t^4}\right)$, on stochastic gradient evaluation (or SAGD iteration) for a fixed system within the stage $t$.

\subsection*{Improving efficiency with same lower and upper level objective}
In the previous section, we deal with the general problem with different lower and upper level objective functions ($h$ and $f$). While the aforementioned methods can be applied to the special case when $f = h$, it misses an opportunity of improving the computational efficiency. 

Specifically, in the aforementioned approach, to control the error of $h(x) - h(x^*)$, we first bound the near-optimal solution $x-x^*$ by bounding the optimization error of the lower objective $f(x)- f(k,x^*)$. When $h=f$, we can directly bound  $h(k,{x}) - h(k,x^*)$. This leads to the following computation of $N^t$ at stage $t$:
  \begin{equation} \label{eq: calculate N_k same}
 \begin{aligned}
      N^t = \min&\left\{ N: \frac{2 \lambda_t \sigma_{G} D }{\sqrt{3N}} + \frac{4M^2 + 4(1+\lambda_t)\sigma_{G}^2}{\mu(N+1)} \right.\\
      &\left.~~~~~~~~+ \frac{4\nu D^2}{N(N+1)} \le \epsilon_t \right\}
 \end{aligned}
 \end{equation}
 Here we replace the right hand side of the inequality $\frac{\mu \epsilon_t^2}{2 L^2}$ with $\epsilon_t$. 
 Consequently, this leads to $N^t = O\left( \frac{\log^2\left(\frac{T}{\alpha}\right)}{\epsilon_t^2}\right)$, which improves the previous order $O\left( \frac{\log^2\left(\frac{T}{\alpha}\right)}{\epsilon_t^4}\right)$
 by $\epsilon_t^2$.

The details of the \textbf{Optimization} step is presented in Algorithm \ref{alg: solution update}.

\begin{algorithm}[t]
   \caption{Optimization}
   \label{alg: solution update}
\begin{algorithmic}[1]
   \STATE {\bfseries Input:} Remaining set $\calR_t$; initial solution $x_{k,N_{k}^{t-1}}$ for $k\in \calR_t$; number of SAGD conducted $N_k^{t-1}$; and tolerance $\epsilon_t$.
   \STATE For each $k\in\calR_t$, compute $N_k^t$ using \eqref{eq: calculate lambda} and \eqref{eq: calculate N_k} if $h(k,\cdot)\neq f(k,\cdot)$ or using \eqref{eq: calculate lambda} and \eqref{eq: calculate N_k same} if $h(k,\cdot)= f(k,\cdot)$.
  \FOR{each $k\in\calR_{t}$} 
    \FOR{$\ell$ from $N_k^{t-1}+1$ to $N^t_k$ }
    \STATE Obtain a stochastic gradient $G(k,x_{k,\ell},\xi_{k,\ell})$ and update $x_{k,\ell}$ using \eqref{eq: SAGD update}.
    \ENDFOR
    \ENDFOR
   \STATE {\bfseries Output:} $x_{k,N_{k}^{t}}$ for all $k\in\mathcal{R}_t$. 
\end{algorithmic}
\end{algorithm}

\subsection{More Efficient Optimization step with Asymptotic Concentration Bound} \label{sec: optimization asymptotic}


In the previous section, we employ SAGD as the optimization tool and its finite-sample error bound to determine the total number of iterations of SAGD. While this approach possesses good statistical guarantee that both the (intermediate) confidence level and the (intermediate) tolerance level are satisfied, it can be computationally inefficient, meaning that the required iteration number $N_k^t$ computed from the concentration bound in Lemma \ref{lem: SAGD concentration} is unnecessarily large. 

In this section we aim to improve the computational efficiency by developing asymptotically valid concentration bound for the stochastic gradient descent (SGD) that results in better sample complexity.
In particular, we develop two asymptotically valid concentration bounds for cases with different and same objectives in lower-level and upper-level problems, respectively. We will show both of them significantly improve the computational efficiency (i.e., total number of iterations $N_k^t$), by a factor of $\frac{\log(\frac{T}{\alpha})}{\epsilon^2_t}$ for the case of different objectives and by a factor of $\frac{1}{\epsilon_t}$ for the case of same objective, compared with the non-asymptotic approach in the previous section.
\subsection*{Different Objectives in Lower and Upper Levels}
In this section, we first consider the case of different objectives in lower-level and upper-level problems. As mentioned above, we switch to the simpler stochastic gradient descent (SGD) method. We use the same notation $N^t_k$ to denote the total number of iterations of running SGD for system $k$ up to stage $t$.
Recall the goal is to control the error of the upper-level objective $h(k,{x}_{N^t_k}) - h(k,x_k^*)$, and thus, we want to characterize the asymptotic normality of this error. 
For this purpose, we make the following mild assumptions of regularity conditions. Let $G(k,x_k,\xi)$ denote a stochastic gradient of $F(k,x_k,\xi)$. We make the following assumptions of regularity conditions on $f$ and $h$.
\begin{assumption} \label{assump: asymptotic} \ 
\begin{enumerate}
    \item $f(k,x_k)$ is twice continuously differentiable in $x_k \in \mathcal{X}_k$.
    \item $h(k,x_k)$ is twice differentiable in $x_k$, $\nabla h(k,x_k^*) \neq 0$ and $\sup_{x_k \in \mathcal{X}_k} \|\nabla^2 h(k,x_k)\|_2 < \infty$.
    \item $\mathbb{E}_{\xi}[G(k,x_k,\xi)] = \nabla f(k,x_k)$ and $\{G(k,x_k,\xi_k)\}_{x_k \in \mathcal{X}_k}$ is uniformly integrable.
    \item There exists $\epsilon>0$, such that $\sup_{x_k \in \mathcal{X}_k} \mathbb{E}_{\xi}[\|G(k,x_k,\xi)\|_2^{2+\epsilon}] < \infty $.
    \item $\operatorname{Cov}(G(k,x_k,\xi))$ is continuous in $x_k$, where $\operatorname{Cov}(v)$ is the covariance matrix of the random vector $v$.
\end{enumerate}
\end{assumption}

Since  the remaining of the section focuses on a fixed system $k$, we drop $k$ from notations for simplicity. SGD updates the solution by the following rule:
\begin{equation} \label{eq: SGD update}
    x_{\ell+1}=\arg \min _{x \in \mathcal{X}} \gamma_{ \ell}<G\left( x_{ \ell}, \xi_{ \ell}\right), x>+\frac{\|x-x_{\ell}\|_2^2}{2},
\end{equation}
where $<\cdot,\cdot>$ is the vector inner product and $\gamma_{\ell}$ is the step size of $\ell$th iteration. We choose a decreasing step size $\gamma_{\ell} = \gamma \ell^{-1}$. Notably, our results can be extended to more general step size such as $\gamma_{\ell} = \gamma \ell^{-\beta}$ for some $\beta \in (1/2,1]$. However, $\beta=1$ yields the fastest convergence rate of the solution given by SGD. As a result, we only focus on the choice of $\beta=1$ in the following.  

The asymptotic normality for $h(x_N) - h(x^*)$ is stated in Theorem \ref{thm:asymptotic normality} 
\begin{theorem}\label{thm:asymptotic normality}
Suppose Assumption \ref{assump: Lipschitz-convex}.1, Assumption \ref{assump: convex}.3 and Assumption \ref{assump: asymptotic} hold. For $\gamma > \frac{1}{2 \mu}$,
\begin{equation}
    N^{\frac{1}{2}} (h(x_{N}) - h(x^*)) \Rightarrow \mathcal{N}(0, \tilde{\sigma}^2),
\end{equation}
where ``$\Rightarrow$" denotes convergence in distribution, $\tilde{\sigma}^2 =  \nabla h(x^*)^\top \Sigma_{\infty} \nabla h(x^*)$, and 
$$\operatorname{vec}(\Sigma_{\infty}) = \gamma^2((\gamma \nabla^2f(x^*)-\frac{I}{2}) \oplus (\gamma \nabla^2 f(x^*)-\frac{I}{2}))^{-1} \operatorname{vec}(\operatorname{Cov}(G(x^*,\xi))).$$
Here $\operatorname{vec}$ is the stack operator and $\oplus$ is the Kronecker sum. Furthermore, 

\centerline{$\Tilde{\sigma}^2 \le \frac{\gamma^2}{2\gamma\mu-1}\|\nabla h(x^*)\|_2^2 \|\nabla^2 f(x^*)\|_2 \|\operatorname{Cov(G(x^*,\xi))\|_2}$.}
\noindent The right hand side is minimized at $\gamma^* = \frac{1}{\mu}$.
\end{theorem}

In Theorem \ref{thm:asymptotic normality}, we derive the weak convergence of $h(x_{N})-h(x^*)$ as well as provide a guidance on how to choose the step size by minimizing an upper bound on the asymptotic variance of the normalized error.
Theorem \ref{thm:asymptotic normality} indicates the error is approximately normally distributed when $N$ is sufficiently large, and hence, we have approximately
$$ \mathbb{P}\left( \left|h(x_{N}) - h(x^*) \right| \ge \epsilon_t)\right) = \mathbb{P}\left( \left|\frac{\sqrt{N}}{\Tilde{\sigma}} \left(h(x_{N}) - h(x^*)\right)\right| \ge \frac{\sqrt{N}}{\Tilde{\sigma}} \epsilon_t)\right) \le \frac{\sqrt{2 } \Tilde{\sigma} }{\sqrt{\pi N} \epsilon_t } \exp\left( -\frac{N \epsilon_t ^2 }{2\Tilde{\sigma}^2}\right).$$
Accordingly, we can choose 
\begin{equation} \label{eq: calculate N_k asymptotic}
    N^t := \min \left\{ N: \frac{\sqrt{2 } \Tilde{\sigma} }{\sqrt{\pi N} \epsilon_t } \exp\left( -\frac{N \epsilon_t ^2 }{2\Tilde{\sigma}^2}\right) \le \alpha_t \right\},
\end{equation}
where $\alpha_t = \frac{\alpha}{2 T K}$, and then conduct $N^t - N^{t-1}$ iterations of SGD in stage $t$.


\textbf{Remark 1: } Using the non-asymptotic bound gives the number of iterations needed for SAGD of the order $N^t = O\left( \frac{\log^2\frac{T}{\alpha}}{\epsilon_t^4}\right)$, whereas using the asymptotic valid bound gives the number of iterations needed for SGD of order $N^t = O\left( \frac{\log\frac{T}{\alpha}}{\epsilon_t^2}\right)$, improving by a factor of the $O\left( \frac{\log\frac{T}{\alpha}}{\epsilon_t^2}\right)$.

\textbf{Remark 2: } In computing $N^t$, one needs to know the value of $\|\nabla h(x^*)\|_2$, $\|\nabla^2 f(x^*)\|_2$ and $\|\operatorname{Cov}(G(x^*,\xi)\|_2$, all of which are unknown since $x^*$ is unknown. In practice, one possible approach is to first run a few SGD iterations to get $\Bar{x}$ as an estimate of the optimal solution $x^*$ and then evaluate these unknown values under solution $\Bar{x}$ using, for example, Monte Carlo simulation. One can also update the estimate as solution $x$ is updated. Another approach is to replace these unknown values with an upper bound if such prior knowledge is available. Notably, this will not affect the order of the sample complexity discussed in the previous remark. Same argument follows for the case where the objective is the same for the lower and upper levels, which is studied in the following paragraph.

\subsection*{Same Objective in Lower and Upper Levels}
When $h$ is different from $f$, we have that $h({x}_{N})$ converges in an asymptotic order of $ N^{-\frac{1}{2}}$. Theorem \ref{thm:asymptotic normality} shows that, when $h=f$, the asymptotic variance $\Tilde{\sigma} = \nabla h(x^*) \Sigma_{\infty} \nabla h(x^*)^\top = \nabla f(x^*) \Sigma_{\infty} \nabla f(x^*)^\top = 0$ since $x^*$ minimizes $f(\cdot)$. This implies, $h({x}_{N}) = f(x_N)$ converges in an order faster than $ N^{-\frac{1}{2}}$. For an intuitive explanation, express 
$h({x}^*)  \approx h({x}_{N}) + \nabla h({x}_{N})^\top ( x^* - {x}_{N} )$. If $\nabla h(x^*) != 0$ (i.e., $h\neq f$), then $ \nabla h({x}_{N}) $ is also bounded away from $0$ when ${x}_{N}$ is close to $x^*$. Consequently, 
the convergence speed of $h({x}_{N})$ is determined by the convergence speed of $ {x}_{N} - x^*$. However, when $\nabla h(x^*) = \nabla f(x^*) = 0$, both $ \nabla h({x}_{N}) $ and ${x}_{N}-x^*$ converge to $0$, which implies the convergence speed of $h(x_N) = f({x}_{N})$ is faster than that of  $ {x}_{N}$. Hence, we need some extra effort to characterize the convergence speed of $f({x}_{N})$. We make the following one additional assumption of regularity condition on $f(\cdot)$.
 
\begin{assumption} \label{assump: asymptotic same objective} \ 
$f(x)$ is three-times continuously differentiable in $x\in\mathcal{X}$.
\end{assumption}

With Assumption \ref{assump: asymptotic same objective}, we develop the weak convergence of $f(x_N) - f(x^*)$ in the following Theorem \ref{thm: asymptotic normality same objective}.
\begin{theorem} [Weak convergence of $f(x_{N})$] \label{thm: asymptotic normality same objective}

Suppose Assumption \ref{assump: Lipschitz-convex}.1, \ref{assump: convex}.3, Assumption \ref{assump: asymptotic} and Assumption \ref{assump: asymptotic same objective}  hold. Then, for $\gamma > \frac{1}{2\mu}$,
    \begin{equation}
       N (f(x_{N}) - f(x^*)) \Rightarrow  Z^\intercal \nabla^2 f(x^*) Z,
    \end{equation}
    where $Z \sim\mathcal{N} \left(0, \Sigma_{\infty} \right)$ and
    $\operatorname{vec}(\Sigma_{\infty})$ is defined in Theorem \ref{thm:asymptotic normality}.
    
\end{theorem}

Theorem \ref{thm: asymptotic normality same objective} implies that, the objective function of the lower level problem converges in an rate of $O(\frac{1}{N})$, which is faster by a factor of $\frac{1}{\sqrt{N}}$ than the convergence rate of the upper-level objective function in the case when the two objective functions are different. Furthermore, the distribution of the normalized error $N(f(x_{N})-f(x^*))$ does not converge to a normal distribution but the distribution of a quadratic form of normal random vectors. However, the non-normal limiting distribution makes determining the number of SGD iterations to reach the desired tolerance and confidence level more complicated as we need (tight) concentration bound for quadratic of normal random vectors.  We summarize the result of how to determine $N$ in the following theorem.

\begin{theorem}\label{thm: quadratic normal concentration}
    Suppose $ N(f(x_{N}) - f(x^*)) \sim Z^\intercal \nabla^2 f(x^*) Z $, where $Z \sim \mathcal{N} \left(0, \Sigma_{\infty} \right)$. Then,  if we set $N^t$ to be
    \begin{equation}\label{eq: calculate N_k asynptotic same}
    \begin{aligned}
         \min\left\{N: N\ge \frac{b}{\epsilon_t} \max\left\{\left[4\log\frac{1}{\alpha_t}+\frac{3d}{2}\right], 2d \right\} \right\},
    \end{aligned}
    \end{equation} we have 
    $$\mathbb{P}\left( f(x_{N^t}) - f(x^*) > \epsilon_t \right) \le \alpha_t,$$
    where recall $d$ is the dimension of $x$, and $b = \frac{\gamma^2}{2\gamma\mu-1} \| \operatorname{Cov}(G(x^*,\xi))\|_2\|\nabla^2 f(x^*)\|_2$. Furthermore, $N^t$ is minimized at $\gamma^* = \frac{1}{\mu} $ and $b^* = \frac{1}{\mu^2}\| \operatorname{Cov}(G(x^*,\xi))\|_2\|\nabla^2 f(x^*)\|_2 $.
\end{theorem}

\textbf{Remark:} Using the non-asymptotic concentration bound gives the number of iterations needed for SAGD of order $N^t = O\left(\frac{\log^2\frac{T}{\alpha}}{\epsilon^2_t}\right)$.
 Theorem \ref{thm: quadratic normal concentration} implies when the objective functions of the lower level and the upper level are the same, the number of SGD iterations needed is of order $N^t = O\left(\frac{\log\frac{T}{\alpha}}{\epsilon_t}\right)$, improving by a factor of $\frac{\log\frac{T}{\alpha}}{\epsilon_t}$.

The  detailed algorithm of the \textbf{Optimization} step with asymptotic approach is presented in Algorithm \ref{alg: solution update asymptotic}.

\begin{algorithm}[t]
   \caption{Optimization (Asymptotic)}
   \label{alg: solution update asymptotic}
\begin{algorithmic}[1]
   \STATE {\bfseries Input:} Remaining set $\calR_t$; initial solution $x_{k,N_{k}^{t-1}}$ for $k\in \calR_t$; number of SAGD conducted $N_k^{t-1}$; and tolerance $\epsilon_t$.
   \STATE For each $k\in\calR_t$, compute $N_k^t$ using \eqref{eq: calculate N_k asymptotic} if $h(k,\cdot)\neq f(k,\cdot)$ or using \eqref{eq: calculate N_k asynptotic same} if $h(k,\cdot)= f(k,\cdot)$.
  \FOR{each $k\in\calR_{t}$} 
    \FOR{$\ell$ from $N_k^{t-1}+1$ to $N^t_k$ }
    \STATE Obtain a stochastic gradient $G(k,x_{k,\ell},\xi_{k,\ell})$ and update $x_{k,\ell}$ using \eqref{eq: SGD update}.
    \ENDFOR
    \ENDFOR
   \STATE {\bfseries Output:} $x_{k,N_{k}^{t}}$ for all $k\in\mathcal{R}_t$. 
\end{algorithmic}
\end{algorithm}

\section{Statistical Validity}
\label{sec:validity}

In this section, we show the statistical validity of the proposed algorithm. Specifically, we show with probability at least $1-\alpha$, Algorithm \ref{alg: RSCOP} outputs an $\epsilon$-optimal system.

To begin with, the following Lemma \ref{lem: optimization validity} guarantees with high probability, the {\bf Optimization} step with the non-asymptotic approach (Algorithm \ref{alg: solution update}) returns an $\epsilon_t$-optimal decision for each system in the remaining set $\calR_t$.

\begin{lemma} \label{lem: optimization validity}
    Suppose Assumption \ref{assump: convex} holds. With probability at least $1-\frac{\alpha}{2}$, the decision variable ${x}_k^t = {x}_{k,N_k^t}$ returned by Algorithm \ref{alg: solution update} satisfies
    $h(k,{x}_k^t) \le h(k,x_k^*) + \epsilon_t$ for $\forall 1\le t\le T, k \in \calR_t $.
\end{lemma}
Notably, although in Section \ref{sec: optimization asymptotic} we claim Algorithm \ref{alg: solution update asymptotic} (the \textbf{Optimization} step with asymptotic approach) leads to better sample complexity, unfortunately it does not always guarantee the required tolerance and confidence level as the concentration bound we develop for Algorithm \ref{alg: solution update asymptotic} is only valid in the asymptotic sense. The efficiency using Algorithm \ref{assump: asymptotic same objective} as the \textbf{Optimization} step will be demonstrated 
numerically in Section \ref{sec:numerical}.

Next, we establish the validity of the \textbf{Pruning} step in Lemma \ref{lem: prunning validity}. Specifically, we aim to determine what can be inferred about the expected performance, $h(k, x_k^t)$, after the \textbf{Pruning} step at stage $t$, based on whether system $k$ remains in $\calR_t$.

\begin{lemma} \label{lem: prunning validity}
    Under Assumption \ref{assump: prunning}, conditioned $\calR_t, x_k^t,\forall k\in \calR_t$, with probability at least $1-\frac{\alpha}{2T}$, 
    \begin{enumerate}
        \item if $k\in\calR_{t+1}$, then $\forall i\neq k \in \calR_{t}$,
    $h(k,x_k^t) - h(i,x_i^t) \le \epsilon'_t$;
    \item otherwise if $k\in\calR_t \backslash\calR_{t+1} $, there exists $i\in\mathcal{R}_{t+1}$,
    $ h(k,x_k^t) - h(i,x_i^t) \ge \epsilon_t$.
    \end{enumerate}
    \end{lemma}
Lemma \ref{lem: prunning validity} indicates that on the one hand, if a system remains in $\calR_{t+1}$ after the \textbf{Pruning} step at stage $t$, we must have, with high probability, its upper-level objective value within $\epsilon'_t$ compared to any other system in $\calR_{t+1}$ under current solutions;  on the other hand, if a system is eliminated by \textbf{Pruning} step at stage $t$, then with high probability its upper-level objective is at least $\epsilon_t$ larger than some other remaining system under current solutions. Notably, this is slightly different from our motivation on examining the feasibility of \eqref{eq: eliminate i} and \eqref{eq: eliminate k}, where $\epsilon_t$ does not show up. In fact, if the true difference on the expected upper level objective between two different systems $h(i,x_i^t)-h(k,x_k^t)$ (or $h(k,x_k^t)-h(i,x_i^t)$ ) lies in the interval $(\epsilon_t,\epsilon'_t)$,  \eqref{eq: eliminate i} (or \eqref{eq: eliminate k}) may be determined either feasible or infeasible. In fact, the difference $\epsilon'_t - \epsilon_t$ represents the acceptance level on the feasibility check itself. A large acceptance level can help reduce the computational effort when some of the constraints are very close to the equality, where a strict feasibility check can be extremely hard. 

With the help of Lemma \ref{lem: optimization validity} and Lemma \ref{lem: prunning validity} guaranteeing the statistical validity of the \textbf{Optimization} and \textbf{Pruning} steps, we can provide the validity result of Algorithm \ref{alg: RSCOP} in the following Theorem.

\begin{theorem} \label{thm: statistical validity}
    Suppose Assumptions \ref{assump: convex} and \ref{assump: prunning} hold. Then, with probability at least $1-\alpha$, Algorithm \ref{alg: RSCOP} outputs an $\epsilon$-optimal system.
\end{theorem}

\section{Numerical Study: Drug selection with optimal dosage}
\label{sec:numerical}
We consider the problem of drug selection with optimal dosage that is adapted from \cite{verweij2020randomized,si2022selecting}. Suppose there are $K=20$ different drugs(treatment) for the same deceases to be selected. For each drug $1\le i\le K$, an optimal dosage $x_i^* \in [0,2]$ needs to be decided to obtain the best expected effect, which is measured by the mean change from baseline in sitting diastolic blood pressure (SiSBP, see \cite{verweij2020randomized}).
The optimal dosage is chosen to minimize SiSBP. The strong convexity of SiSBP in terms of the drug dosage has been examined by empirical evidence (e.g., a quadratic form of SiSBP was identified in \cite{verweij2020randomized}). Since the data in \cite{verweij2020randomized} is not public, we study a synthetic setting. To be specific, for each drug $i$, suppose the SiSBP has the form
$$f(i,x_i) = a_{i,2}x_i^2 + a_{i,1}x_i + a_{i,0},$$
where we set $a_{i,2} = 1+0.1*i,a_{i,1} = -3*a_{i,2}, a_{i,0} = \frac{a^2_{i,1}}{4a_{i,2}} + 0.11*i$ for $i=1,\ldots,K$. That is, the true optimal dosage $x_i^* =1.5$ and the corresponding optimal SiSBP is $f(i,x_i^*) = 0.11*i$.
Furthermore, we consider the stochastic counter part
$$F(i,x_i,\xi) = (a_{i,2} + \xi_{2}) x_i^2 + (a_{i,1}+\xi_1) x_i + (a_{i,0}+\xi_0),$$
where $\xi_r \sim \operatorname{Uniform}(-0.5,0.5), r= 0,1,2$. The intuition of the stochastic function $F$ is to model the potential discrepancies between individuals within a population, which is represented by the stochastic (perturbed) parameters $a_{i,r},r=0,1,2$. When testing drug effectiveness with a certain dosage, one has no access to exact evaluation of $f$ but stochastic samples of  $F$. 

We consider two different scenarios with the same or different upper and lower level objective functions. In the first scenario, one tries to find the drug with the best expected effect under its optimal dosage, i.e., $\min_{i} f(i,x_i^*), ~x_i^*=\arg\min_{x_i \in [0,2]} f(i,x_i)$. In the second scenario, while for each drug the optimal dosage needs to be decided to obtain the minimal SiSBP, one also considers the cost (or reward) of providing the drug (treatment). As a result, the problem becomes a bilevel optimization with the upper level choosing the best drug that minimizes a weighted summation of the cost and effectiveness and the lower level of purely minimizing the SiSBP for each drug. Mathematically, the upper level objective is defined as 
$h(i,x_i) = c_i x_i + \omega_i f(i,x_i)$ and the stochastic counter part $H(i,x_i,\xi) = c_i x_i + \omega_i F(i,x_i,\xi)$. In practice, the weight $\omega_i$ depends on the preference over cost and effect of drugs. In this numerical study, we set $\omega_i = 1$ and $c_i = 1$. 

For implementation, $\epsilon = 0.1$ for both different objectives and for the same objective, both of which insures that only the true optimal system satisfies $\epsilon$-near optimality (The performance difference between the best and second best system is $0.11$). The confidence level is set to $0.9$. We run $500$ macro-replications to obtain the probability of selecting an $\epsilon$-optimal system (Probability), average number of function evaluations (Function), and average number of gradient evaluations (Gradient) along with their $95\%$ confidence interval. 


For comparison baselines, to our best knowledge, only \cite{si2022selecting} and this paper consider the generalized R\&S with an optimized continuous decision variable. Moreover, the method by \cite{si2022selecting} can only be applied when a computing budget (total number of SGD iterations) is given ahead while our method determines this computing budget with the given tolerance $\epsilon$ as well as the significance level $\alpha$. For these reasons, it is difficult to compare our method with \cite{si2022selecting}. Therefore, we only test our method with either exact or asymptotic concentration bound for optimization and with different choices of the total number of stages $T = 1,3,4,5$. For the sequence of tolerance values, we set $\epsilon_t = \frac{2}{5} 2^{T-t}\epsilon$ and $\epsilon_t' = \frac{3}{5}2^{T-t} \epsilon, t=1,\ldots,T$.  Note that $T=1$ corresponds to the direct approach in Section \ref{sec: direct approach}. For more implementation details, for the exact approach, we set $\mu_i = 2a_{i,2}$, $\nu_i = \mu_i$, $M_i = 0$, $L_i = 2a_{i,2}+a_{i,1}$, $D_i = 2$, $\sigma_{G,i}^2 = \frac{1}{3}$; for the asymptotic approach, we compute a upper bound on $\tilde{\sigma}_i^2 \le  \frac{c_i D_i}{\mu_i^2}$.

 \renewcommand{\arraystretch}{1.5}

\begin{table}[h]
    \centering
    \tiny
    \begin{tabular}
{|c|c|c|c|c|c|c|} \hline
Method& & T=1& T=2& T=3& T=4& T=5 \\ \hline
\multirow{3}{*}{Exact approach}&Probability& $1.0  (\pm 0.0) $& $1.0  (\pm 0.0) $& $1.0  (\pm 0.0) $& $1.0  (\pm 0.0) $& $1.0  (\pm 0.0) $ \\
&Gradient & $9.02 ( \pm0.0)\times 10^6$& $1.94 ( \pm0.03)\times 10^6$& $1.58 ( \pm0.06)\times 10^6$& $1.69 ( \pm0.06)\times 10^6$& $1.81 ( \pm0.04)\times 10^6$ \\
&Function& $1.78 ( \pm0.46)\times 10^5$& $2.3 ( \pm0.6)\times 10^5$& $2.08 ( \pm0.39)\times 10^5$& $2.28 ( \pm0.42)\times 10^5$& $2.55 ( \pm0.41)\times 10^5$ \\ \hline
\multirow{3}{*}{Asymptotic approach}&Probability& $0.914 (\pm 0.02) $& $0.946 (\pm 0.02) $& $0.946 (\pm 0.02) $& $0.906 (\pm 0.03) $& $0.93 (\pm 0.02) $ \\
&Gradient& $1.21 ( \pm0.0)\times 10^3$& $7.04 ( \pm1.1)\times 10^2$& $5.91 ( \pm1.86)\times 10^2$& $5.67 ( \pm2.07)\times 10^2$& $5.94 ( \pm2.18)\times 10^2$ \\
&Function & $2.51 ( \pm1.99)\times 10^5$& $1.92 ( \pm1.1)\times 10^5$& $1.68 ( \pm0.99)\times 10^5$& $1.69 ( \pm1.14)\times 10^5$& $1.95 ( \pm1.25)\times 10^5$ \\ \hline
\end{tabular}
    \caption{Different Objective}
    \label{tab: different objective}
\end{table}

As Table \ref{tab: different objective} indicates,  in terms of the empirical Probability, Algorithm \ref{alg: RSCOP} with both exact approach or asymptotic approach achieves the target confidence level $0.9$, verifying its statistical validity. As for the number of gradient evaluations, i.e., the total number of SAGD (SGD) iterations run by the algorithm, the exact approach requires more than $10^6$ number of stochastic gradient evaluations, whereas the asymptotic approach requires only $10^2 \sim 10^3$ number of stochastic gradient evaluations, which aligns with the result of less sample complexity by the asymptotic approach. Furthermore, for the empirical Probability, asymptotic approach also reaches a empirical Probability closer to the target value, whereas the  exact approach always yields an empirical Probability equal to $1$, indicating the conservativeness of exact approach which leads to unnecessary optimization effort. The difference in the total number of function evaluations between the exact approach and the asymptotic approach is much less significant as both use the same pruning sub-algorithm. 

For different choices of $T$, the number of tolerance level, we see that as $T$ increases, both the total number of gradient evaluations and the total number of function evaluations first decreases and then increases, indicating a trade-off in the choice of $T$.  In this example, the total number of gradient evaluations for both the exact and asymptotic approach reaches the minimal at $T=4$, and the total number of function evaluations for both the exact and asymptotic approach reaches the minimal at $T=3$. Intuitively speaking, a larger $T$ helps identify some sub-optimal systems with less computational effort due to a larger tolorance value, but may result in more computational effort as one needs to divide the confidence level $\alpha$ by a factor of $T$ when deciding the total number of SAGD (SGD) required to reach the stage-wise desired tolerance and confidence level.

\begin{table}
\tiny
    \centering
    \begin{tabular}{|c|c|c|c|c|c|c|} \hline
Method& &T=1&T=2&T=3&T=4&T=5 \\ \hline
\multirow{3}{*}{Exact approach}&Probability& $1.0  (\pm 0.0) $& $1.0  (\pm 0.0) $& $1.0  (\pm 0.0) $& $1.0  (\pm 0.0) $& $1.0  (\pm 0.0) $ \\
&Gradient & $2.34 ( \pm0.0)\times 10^3$& $1.23 ( \pm0.0)\times 10^3$& $9.33 ( \pm0.06)\times 10^2$& $9.17 ( \pm0.13)\times 10^2$& $9.68 ( \pm0.14)\times 10^2$ \\
&Function & $1.78 ( \pm0.45)\times 10^5$& $1.76 ( \pm0.41)\times 10^5$& $1.68 ( \pm0.29)\times 10^5$& $1.89 ( \pm0.3)\times 10^5$& $2.11 ( \pm0.35)\times 10^5$ \\ \hline
\multirow{3}{*}{Asymptotic approach}&Probability& $0.992 (\pm 0.01) $& $0.998 (\pm 0.0) $& $0.998 (\pm 0.0) $& $1.0 (\pm 0.0) $& $0.996 (\pm 0.01) $ \\
&Gradient & $2.54 ( \pm0.0)\times 10^2$& $1.91 ( \pm0.04)\times 10^2$& $1.59 ( \pm0.05)\times 10^2$& $1.54 ( \pm0.05)\times 10^2$& $1.56 ( \pm0.06)\times 10^2$ \\
&Function & $1.78 ( \pm0.88)\times 10^5$& $1.76 ( \pm0.57)\times 10^5$& $1.68 ( \pm0.54)\times 10^5$& $1.89 ( \pm0.43)\times 10^5$& $2.11 ( \pm0.56)\times 10^5$ \\ \hline

    \end{tabular}
    \caption{Same Objective}
    \label{tab: same objective}
\end{table}

Table \ref{tab: same objective} indicates a similar result as Table \ref{tab: different objective}. First, both exact and asymptotic approach achieves the target confidence level, showing their statistical validity. Second, the asymptotic approach is more efficient than the exact approach in the sense that the total number of gradient evaluations by the exact approach is more than that by the asymptotic approach and that asymptotic approach achieves a empirical Probability closer to the target value. Nonetheless, the asymptotic approach also reaches a empirical Probability much higher than the target confidence level $0.9$. Compared with the scenario with different objective, the augment conservativeness in this scenario with the same objective comes from the relaxation in the asymptotic concentration bound in Theorem \ref{thm: quadratic normal concentration}, where we need to develop concentration bound for a quadratic of normal random vector. Unlike normal distribution, where tight concentration bound is available, here we are only able to develop a relatively loose concentration bound (compared to the normal case). Despite the less significant improvement, the asymptotic approach still outperforms the exact approach in terms of the computational efficiency.


\section{Conclusion}
\label{sec:conclusion}
In this paper, we propose a fixed confidence fixed tolerance formulation for solving a stochastic bi-level optimization with categorical upper-level decision variable and continuous lower-level decision variables.    
To solve this problem, we introduce a multi-stage pruning-optimization framework.  The proposed framework alternates between optimizing the decision variable for each system and comparing the performances of different systems. This framework is designed to enhance computational efficiency by identifying and pruning inferior systems with large tolerance early on, thereby avoiding redundant computational efforts associated with optimizing these systems. In optimizing the lower-level decision variable, we propose an asymptotic approach by developing the asymptotic convergence for the objective function under solutions solved using stochastic gradient descent. Compared with a non-asymptotic approach utilizing the existing non-asymptotic concentration bound, our proposed asymptotic approach significantly improve the sample efficiency. We also demonstrate the effectiveness of the framework both theoretically and numerically.


\bibliographystyle{plain}
\bibliography{ref}
\appendix
\section{Stochastic Accelerated Gradient Descent (SAGD)} \label{ECsec: SAGD}
Let $G(k,x_k,\cdot)$ denote the stochastic gradient of $f(k,x_k)$. 
We make the following additional assumption of regularity conditions on $f(k,\cdot)$, which is required by SAGD as in \cite{Lan:book2020}.
\begin{assumption}\label{assump: convex}
For each fixed $k\in\mathcal{K}$,
\begin{enumerate} 
\item $f(k,\cdot)$ satisfies
$$
f(k,y)-f(k,x)-\left\langle \nabla f(k,x), y-x\right\rangle \leq \frac{\nu_k}{2}\|y-x\|_2^2+M_k\|y-x\|_2, \quad \forall x, y \in \mathcal{X}_k,
$$
where $\|\cdot\|_2$ denotes the $l_2$ norm.
    \item  $F(k,x_k,\xi)$ is convex in $x_k$ for every system $k$ and almost every $\xi$.
    \item $\mathbb{E}_{\xi}[G(k,x_k,\xi)] = \nabla f(k,x_k)$.
    \item There exists $\sigma_{G,k}>0$, such that $G(k,x_k,\cdot)$ is $\sigma_{G,k}$-sub-Gaussian random variable. In particular, a random variable $X$ is $\sigma$-sub-Gaussian if 
    $$\mathbb{E}[e^{\lambda^2 X^2} ]\le e^{\lambda^2\sigma^2}, \quad \forall |\lambda| \le \frac{1}{\sigma}.$$
    \item There exists $D_k > 0$, $D^2_k = \sup_{x,y\in \mathcal{X}_k} \|x-y\|_2^2 < \infty$.
\end{enumerate}    
\end{assumption}

On top of Assumption \ref{assump: Lipschitz-convex}.1, Assumption \ref{assump: convex}.1 covers a wide range of strongly convex function. For instance, when $f(k,\cdot)$ is $M_k$-Lipschitz continuous, then we have the second inequality holds with $\nu_k$ = 0. Furthermore, if $\nabla f(k,\cdot)$ is  Lipschitz continuous, then the second inequality holds with $M_k = 0$. Other assumptions are mild regularity conditions. We refer to \cite{Lan:book2020} for detailed discussions on these assumptions.

In addition, since SAGD is applied for each fixed system $k$, in the remaining of the section, we suppres the notation of $k$ in both subscripts and function input (i.e., $f(k,\cdot)$  is simplified as $f(\cdot)$). SAGD updates the solution with the following rule at each iteration $l$:

\textbf{Stochastic Accelerated Gradient Descent (SAGD):}
\begin{equation}\label{eq: SAGD update}
    \begin{aligned}
\underline{x}_{\ell}=&\left(1-q'_{\ell}\right) {x}_{\ell-1}+q'_{\ell} \Bar{x}_{\ell-1}, \\
 \Bar{x}_{\ell}=&\arg \min _{x \in \mathcal{X}}\left\{\gamma_{\ell}\left[\left\langle G\left(\underline{x}_{\ell}, \xi_{\ell}\right), x\right\rangle + \right.\right.\\
 &\left.\left.\mu \frac{\|\underline{x}_{\ell} - x\|_2^2}{2}\right]+\frac{\|\Bar{x}_{\ell-1}- x\|_2^2}{2}\right\}, \\
 {x}_{\ell}=&\left(1-q_{\ell}\right) {x}_{\ell-1}+q_{\ell} \Bar{x}_{\ell},
\end{aligned}
\end{equation}
where $q_{\ell} = \frac{2}{\ell+1}$, $\frac{1}{\gamma_{\ell}} = \frac{\mu(\ell-1)}{2} + \frac{2\nu}{\ell}$ and $q'_{\ell} = \frac{q_\ell}{q_\ell + (1-q_\ell)(1+\mu \gamma_{\ell})}$. Here $\underline{x}_{\ell},\Bar{x}_{\ell}$ are intermediate solutions in the $\ell$th iteration of SAGD and $\{x_{\ell}\}_{\ell=1}^\infty$ is the sequence of solutions generated by SAGD. 

\section{Proof of Theorem \ref{thm:asymptotic normality}}
\textit{Proof.}
We fix a system $k$ and drop the notation on $k$ for simplicity, i.e., $h(x_N)$ is short for $h(k,x_{k,N_k}) $.  Under Assumption \ref{assump: asymptotic}.2, we have 
      $$N^{\frac{1}{2}}\left( h(x_N) - h(x^*) \right) 
      = \underbrace{N^{\frac{1}{2}}\nabla h(x^*) ^\top(x_N-x^*)}_{Z_{N,1}}
        +  \underbrace{N^{\frac{1}{2}}(x_N-x^*)^\top C(x_N) (x_N-x^*) }_{Z_{N,2}},$$
    where $C(x_N)$ is some matrix satisfying $\|C(x_N)\|_2 \le \sup_x \|\nabla^2 h(x) \|_2 < \infty$. 
    \begin{itemize}
        \item  $Z_{N,1}$: 
        Let $U_N = N^{\frac{1}{2}}(x_N-x^*)$ and $t_N =\sum_{n=1}^N n^{-1}$.
        Under Assumption \ref{assump: convex}.1 and Assumption \ref{assump: asymptotic}, by Theorem 10.2.1 in \cite{Kushner2003}, we have 
        $U_N$ converge weakly to $U(t_N)$, where $U(t)$ is the solution to the following stochastic differential equations 
        $$ 
           \mathrm{d}U = -(\gamma \nabla^2f(x^*) - \frac{1}{2}I)  U\mathrm{d}t + \operatorname{Cov}^{\frac{1}{2}}(G(x^*,\xi)) \mathrm{d} B_t$$
        Here $B_t$ is the standard Brownian motion with an identity covariance matrix and $\operatorname{Cov}^{\frac{1}{2}}(G(x^*,\xi))$ is a positive definite matrix satisfying $\left(\operatorname{Cov}^{\frac{1}{2}}(G(x^*,\xi))\right)^2 = \operatorname{Cov}(G(x^*,\xi))$. The stochastic process satisfying the above stochastic differential equation belongs to the family of  Ornstein-Uhlenbeck (OU) process (e.g., see \cite{meucci2009review}). For a fixed $t$, $U(t)$ follows a multivariate normal distribution with mean $0$ and covariance matrix $\Sigma_t$ that converges to 
  $\Sigma_\infty$, where $$ 
           \Sigma_\infty = \int_0^{\infty} e^{ \lpar -\gamma \nabla^2f(x^*) +\frac{1}{2}I\rpar t}  \operatorname{Cov}(G(x^*,\xi)) e^{(-\gamma \nabla^2 f(x^*) +\frac{1}{2}I) t} d t
        $$
$\Sigma_\infty$ also has the following explicit expression. 
$$\operatorname{vec}(\Sigma_{\infty}) = \gamma^2((\gamma \nabla^2f(x^*)-\frac{I}{2}) \oplus (\gamma \nabla^2 f(x^*)-\frac{I}{2}))^{-1} \operatorname{vec}(\operatorname{Cov}(G(x^*,\xi))).$$
        This implies $N^{\frac{1}{2}}(x_N - x^*) \Rightarrow \mathcal{N}(0,\Sigma_\infty)$. And hence, $Z_{N,1} \Rightarrow \mathcal{N}(0, \nabla h(x^*)^\top \Sigma_\infty \nabla h(x^*))$.
    \item $Z_{N,2}$:
    Notice $|Z_{N,2}| \le \gamma N^{\frac{1}{2}}\|C(x_n)\|_2 \|x_n-x^*\|^2 \le\gamma \sup_x\|\nabla^2 h(x)\|_2 N^{\frac{1}{2}} \|x_n-x^*\|_2^2$. It suffices to prove $N^{\frac{1}{2}} \|x_n-x^*\|_2^2 \Rightarrow 0 $. To see this, $\forall \varepsilon > 0$, 
    \begin{align*}
        &\mathbb{P} \lpar N^{\frac{1}{2}} \|x_n-x^*\|_2^2 > \varepsilon  \rpar \\
        = & \mathbb{P} \lpar (N^{\frac{1}{2}} \|x_n-x^*\|_2)^2 > \varepsilon N^{\frac{1}{2}}  \rpar\\
        = & \mathbb{P} \lpar N^{\frac{1}{2}} \|x_n-x^*\|_2 > \sqrt{\varepsilon N^{\frac{1}{2}} } \rpar \rightarrow 0,
    \end{align*}
since $N^{\frac{1}{2}}(x_n - x^*)$ converge weakly to $\mathcal{N}(0, \nabla h(x^*)^\top \Sigma_\infty \nabla h(x^*))$, and $\sqrt{\varepsilon N^{\frac{1}{2}} } \rightarrow \infty$. 
Hence, $Z_{N,2} \Rightarrow 0$. 
    \end{itemize}
    Together with the convergence of $Z_{N,1}$ and $Z_{N,2}$, We complete the proof.
\hfill $\blacksquare$

\section{Proof of Theorem \ref{thm: asymptotic normality same objective}}

\textit{Proof.}
Again we fix a system $k$ and drop the notation on $k$ for simplicity.  Under Assumption \ref{assump: asymptotic same objective}, we have
$$      N\left( f(x_N) - f(x^*) \right)=   \underbrace{N(x_N-x^*)^\top \nabla^2f(x^*) (x_N-x^*) }_{Z_{N,1}} + \underbrace{O\lpar N\|x_N - x^*\|_2^3 \rpar}_{Z_{N,2}}.
$$   Since $\mathcal{X}$ is bounded, by the differential continuity of $\nabla^2 f(x)$, we have there exists $0<C<\infty$, $Z_{N,2} \le C\|x_N-x^*\|_2^3$.
    \begin{itemize}
        \item  $Z_{N,1}$: 
        Similar to proof of Theorem \ref{thm:asymptotic normality},
        we have $N^{\frac{1}{2}}(x_N - x^*) \Rightarrow \mathcal{N}(0,\Sigma_\infty)$. And hence, $Z_{N,1} \Rightarrow Z^\top \nabla^2 f(x^*) Z$, where $Z \sim \mathcal{N}(0,\Sigma_\infty)$.
    \item $Z_{N,2}$:
    Also following a similar proof as that of Theorem \ref{thm:asymptotic normality}. $\forall \varepsilon > 0$, 
    \begin{align*}
        &\mathbb{P} \lpar N \|x_n-x^*\|_2^3 > \varepsilon  \rpar \\
        = & \mathbb{P} \lpar (N^{\frac{1}{2}} \|x_n-x^*\|_2)^3 > \varepsilon N^{\frac{1}{2}}  \rpar\\
        = & \mathbb{P} \lpar N^{\frac{1}{2}} \|x_n-x^*\|_2 > \left(\varepsilon N^{\frac{1}{2}} \right)^{\frac{1}{3}} 
        \rpar \rightarrow 0 ~~\text{ as }N\rightarrow \infty.
    \end{align*}
The convergence holds since $N^{\frac{1}{2}}(x_n - x^*)$ converge weakly to $\mathcal{N}(0, \nabla h(x^*)^\top \Sigma_\infty \nabla h(x^*))$, and $\lpar\varepsilon N^{\frac{1}{2}} \rpar^{\frac{1}{3}} \rightarrow \infty$. 
Hence, $Z_{N,2} \Rightarrow 0$. 
    \end{itemize}
Together with the convergence of $Z_{N,1}$ and $Z_{N,2}$, We complete the proof.
\hfill $\blacksquare$

\section{Proof of Theorem \ref{thm: quadratic normal concentration}}
To obtain an asymptotic concentration bound, we first introduce the following result (see, e.g.,  \cite{gallagher2019improved} and \cite{christ2017ancestral}).
\begin{lemma}[Theorem 1 in \cite{gallagher2019improved}] \label{lem: quadratic normal}
Let $X \sim \mathcal{N}(0,I)$ be a standard multi-dimensional normal distribution and $Q = X^\top \Sigma X$, where $\Sigma$ is a positive definite matrix. Let $v = 4\|\Sigma\|_F^2$ and $b = \|\Sigma\|_2$, where $\|\cdot\|_F$ is the Frobenius norm. Then, for $q > \mathbb{E}[Q] + \frac{v}{4b},$
\begin{equation} \label{eq: concentration quadratic gaussian}
\mathbb{P}(Q>q) \le \exp\left(\frac{1}{2}\frac{v}{(4b)^2}\right) \exp\left( -\frac{q-\mathbb{E}[Q]}{4b}\right).
\end{equation}
\end{lemma}
\textit{Proof.}
The goal is to obtain concentration bound for 
$f(k,x_{k,N_k}) - f(k,x_k^*) \ge \epsilon$ for a fixed $k$, $N_k$ and $\epsilon$. Here we drop the notation of stage counter $t$ for simplicity. With the approximation $N_k\lpar f(k,x_{k,N_k}) - f(k,x_k^*) \rpar \sim Z^\top \nabla f(x_k^*) Z$, $Z \sim \mathcal{N}(0,\Sigma_{k,\infty})$,
\begin{align*}
   & \{f(k,x_{k,N_k}) - f(k,x_k^*) \ge \epsilon \} \\
   = &  \{N_k \lpar f(k,x_{k,N_k}) - f(k,x_k^*) \rpar \ge N_k \epsilon \}\\
   = &\{ X^\top \Sigma X \ge N_k \epsilon \}.
\end{align*}
In the last equality, $X$ is a standard multi-variate normal random variable with dimension $d_k$ and $\Sigma = \Sigma^{\frac{1}{2}}_{k,\infty}\nabla^2 f(k,x_k^*) \Sigma^{\frac{1}{2}}_{k,\infty} $. Then, we want to use Lemma \ref{lem: quadratic normal} with this specific $\Sigma$. Nonetheless, we need to bound the parameters $b$ and $v$ in Lemma \ref{lem: quadratic normal} in order to obtain computable concentration bound.

First, notice $\|A\|_F^2 \le d \|A\|_2^2$ for any $A\in\mathbb{R}^{d\times d}$. Since $x_k \in \mathbb{R}^{d_k}$, we have 
$ \frac{v}{4b^2} \le d. $ The first term in \eqref{eq: concentration quadratic gaussian} can be bounded by $\exp(\frac{d}{8})$. Furthermore, since $\Sigma_{k,\infty}$ is symmetric, $\|\Sigma\|_2 \le \|\Sigma^{\frac{1}{2}}_{k,\infty}\|_2^2 \|\nabla^2 f(k,x_k^*)\|_2 =  \|\Sigma^{\frac{1}{2}}_{k,\infty}\|_2 \|\nabla^2 f(k,x_k^*)\|_2$. 
To bound $\|\Sigma_{k,\infty}\|_2$, notice $\Sigma_{k,\infty}$ has the following alternative expression:
$$ \Sigma_{k,\infty} = \gamma_k^2\int_{t=0}^\infty \mathrm{e} ^{(-\gamma_k \nabla^2 f(k,x_k^*) + I/2)t} \operatorname{Cov}(G(k,x_k^*,\xi)) \mathrm{e} ^{(-\gamma_k \nabla^2 f(k,x_k^*) + I/2)t} \mathrm{d}t.
$$
Since $$\|\Sigma_{k,\infty}\|_2 \le \gamma_k^2 \int _{t=1}^\infty \|\mathrm{e} ^{(-\gamma_k \nabla^2 f(k,x_k^*) + I/2)t} \operatorname{Cov}(G(k,x_k^*,\xi)) \mathrm{e} ^{(-\gamma_k \nabla^2 f(k,x_k^*) + I/2)t}\|_2 \mathrm{d}t,$$ we can bound
 $ \|\mathrm{e} ^{(-\gamma_k \nabla^2 f(k,x_k^*) + I/2)t} \operatorname{Cov}(G(k,x_k^*,\xi)) \mathrm{e} ^{(-\gamma_k \nabla^2 f(k,x_k^*) + I/2)t} \|_2$
for each $t$. For
$\|\mathrm{e} ^{(-\gamma_k \nabla^2 f(k,x_k^*) + I/2)t}\|_2$, notably, $\nabla^2 f(k,x_k^*) \succcurlyeq \mu_k I$, which means $-\gamma_k \nabla^2 f(k,x_k^*) + I/2\le (\frac{1}{2} -\gamma_k \mu_k) I $. Hence, $\|\mathrm{e} ^{(-\gamma_k \nabla^2 f(k,x_k^*) + I/2)t}\|_2 \le \mathrm{e}^{(\frac{1}{2}-\gamma_k\mu_k)t} $. Then,
$$  \|\mathrm{e} ^{(-\gamma_k \nabla^2 f(k,x_k^*) + I/2)t} \operatorname{Cov}(G(k,x_k^*,\xi)) \mathrm{e} ^{(-\gamma_k \nabla^2 f(k,x_k^*) + I/2)t} \|_2 \le\mathrm{e}^{(1-2\gamma_k\mu_k)t}\| \operatorname{Cov}(G(k,x_k^*,\xi))\|_2. $$
And we have 
$$\|\Sigma_{k,\infty}\|_2 \le \frac{\gamma_k^2}{2\gamma_k\mu_k-1} \| \operatorname{Cov}(G(k,x_k^*,\xi))\|_2$$
Then, $b \le  \frac{\gamma_k^2}{2\gamma_k\mu_k-1} \| \operatorname{Cov}(G(k,x_k^*,\xi))\|_2\|\nabla^2 f(k,x_k^*)\|_2.$ And we define $b_k := \frac{\gamma_k^2}{2\gamma_k\mu_k-1} \| \operatorname{Cov}(G(k,x_k^*,\xi))\|_2\|\nabla^2 f(k,x_k^*)\|_2$. 

In addition, Lemma \ref{lem: quadratic normal} can be used only if $N_k \ge \mathbb{E}[X^\top \Sigma X] + \frac{v}{4b}$. Notice $ = \operatorname{tr}(\Sigma) + \frac{v}{4b} \le d_k \|\Sigma\|_2 + d_kb = 2d_kb\le 2d_kb_k$. Then we can set $N_k \epsilon\ge 2d_k b_k$.

Moreover, to obtain the desired nominal error $\alpha$, we then let 
$$\exp{\frac{d}{8}} \exp{-\frac{N_k\epsilon - \mathbb{E}[X^\top \Sigma X]}{4 b}}\le \exp{\frac{d_k}{8}} \exp{-\frac{N_k\epsilon - d_kb_k}{4 b_k}} \le \alpha.$$
This can be implied by
$$ N_k \ge \frac{1}{\epsilon} \left[4b_k \lpar\log\frac{1}{\alpha} + \frac{d_k}{8} \rpar + d_k b_k\right].$$
Then, by setting $N_k \ge \max\left \{\left[4b_k \lpar\log\frac{1}{\alpha} + \frac{d_k}{8} \rpar + d_k b_k\right], \frac{2d_kb_k}{\epsilon}\right\} $, we have with probability at least $1-\alpha$, $f(k,x_{k,N_k}) - f(k,x_k^*) \le \epsilon$. Substitute $\epsilon$ and $\alpha$ with $\epsilon_t$ and $\alpha_t$, respectively, we get the desired result.
\hfill $\blacksquare$

\section{Proof of Lemma \ref{lem: optimization validity}}
\textit{Proof.}
    First consider the case where $h\neq f$. By Lemma \ref{lem: SAGD concentration}, for any $k\in K$, with $\alpha_t = \frac{\alpha}{2T}$, $\lambda_t$ defined in \eqref{eq: calculate lambda}, $N_k^t$ defined in \eqref{eq: calculate N_k}, we have 
    \centerline{$\mathbb{P}\left( f(k,x_{k,N_k^t}) - f(k,x_k^*) > \frac{\mu_k \epsilon^2_t }{2 L^2_k} \right) \le \frac{\alpha}{2TK}.$}
    With a union bound on $t=1,\ldots,T$ and $k\in\calK$, we have with probability $1-\frac{\alpha}{2}$, $\forall t, k$,
    \centerline{$|h(k,x_{k,N_k^t}) - h(k,x_k^*)| \le L_k\sqrt{\frac{2(f(k,x_{k,N^t_k}) - f(k,x_k^*))}{\mu_k}}\le \epsilon_t $.}
    Let $A$ denote the event
    
    \centerline{$A = \left\{
    |h(k,x_{k,N_k^t} - h(k,x_k^*)|\le \epsilon_t, \forall 1\le t\le T,k\in \calK
    \right\}$}
    Then $A$ holds true with probability at least $1-\frac{\alpha}{2}$. Under any trajectory in $A$, we also obtain for any $t,k\in\calR_t$,
    $|h(k,x_{k,N_{N_k^t}}) - h(k,x_k^*)|\le \epsilon_t$. 

    For the case $h=f$, with $N_k^t$ defined by \eqref{eq: calculate N_k same}, same argument follows.
\hfill $\blacksquare$

\section{Proof of Lemma \ref{lem: prunning validity}}
\textit{Proof.} 
Consider the following $2$ inequalities with acceptance level $\tau_t$:
\begin{equation}
    h(i,{x}_i^t) - h(k,{x}_k^t) \le q_t. \label{eq: eliminate i with al}
\end{equation}
\begin{equation}
    h(i,{x}_i^t) - h(k,{x}_k^t) \le -q_t. \label{eq: eliminate k with al}
\end{equation}

The feasibility check in Algorithm \ref{alg: pruning efficient} is done for \eqref{eq: eliminate i with al} and \eqref{eq: eliminate k with al}. 

For a pair $i<k\in\calR_t$, let $A_{ik}^t$ represent the event of \textbf{correct detection} for \eqref{eq: eliminate i with al} at stage $t$. We define $A_{ik}^t$ as follows:  if $y_i^t-y_k^t \ge \epsilon'_t$, the algorithm declares \eqref{eq: eliminate i with al} infeasible; else if $y_i^t-y_k^t \le \epsilon_t$, the algorithm declares \eqref{eq: eliminate i with al} feasible; otherwise, the algorithm declares \eqref{eq: eliminate i with al} either feasible or infeasible. Similarly, we let $B_{ik}^t$ denote the event of \textbf{correct detection} made for \eqref{eq: eliminate k with al} at stage $t$. Under $B_{ik}^t$, if $y_i^t-y_k^t { \le} -\epsilon'_t$, the algorithm declares \eqref{eq: eliminate k with al} feasible; else if $y_i^t-y_k^t { \ge} -\epsilon_t$, the algorithm declares \eqref{eq: eliminate k with al} infeasible; otherwise, the algorithm declares \eqref{eq: eliminate k with al} either feasible or infeasible. Furthermore, we let $\operatorname{CD}^t:= \cap_{i<k\in\calR_t} \left(A_{ik}^t\cap B_{ik}^t\right)$ denote the event of \textbf{correct detection} at stage $t$. 
If $|\calR_t|=1$, we do not need to do the feasibility check and hence 
$ \mathbb{P}\left( \operatorname{CD}^t \mid \calR_t,\mathbf{\Bar{x}}^t\right) = 1 \ge 1 - \frac{\alpha}{2T}.$
 If $|\calR_t| \ge 2$, 
 conditioned on $\calR_t$ and $\mathbf{{x}}^t$,since the random samples of function evaluations are i.i.d. conditioned on $\calR_t$ and $x^t$. By Theorem 2 in \cite{zhou2022finding} with $2$ thresholds for each constraint,  

$\mathbb{P}\left( \operatorname{CD}^t \mid \calR_t,\mathbf{\Bar{x}}^t\right) \ge 1 - \frac{\alpha}{2T}.$ Notably, In Algorithm \ref{alg: pruning efficient}, if $k \in \calR_{t+1}$, then it must satisfy $h(k,x_k^t) - h(i,x_i^t) \le q_t$ is determined feasible for $i>k,i\in \calR_t$ and $h(i,x_i^t) - h(k,x_k^t) \le -q_t$ is determined infeasible for $i<k,i\in\calR_t$. The two cases together indicate, $h(k,x_k^t) \le h(i,x_i^t) + q_t+\tau_t = \epsilon'_t, \forall i\in\calR_t$.  We prove 1.

For 2., since $k\not\in \calR_{t+1}$, there must exist $i\in\calR_t$, either $i>k,\calR_t$ and $h(k,x_k^t) - h(i,x_i^t) \le q_t$ is determined infeasible or $i<k,i\in\calR_t$ and  $h(i,x_i^t) - h(k,x_k^t) \le -q_t$ is determined feasible. Both cases alone indicate there exists $i\in\calR_t$, $h(k,x_k^t) \ge h(i,x_i^t) + q_t - \tau_t = \epsilon_t$. The proof is complete.
\hfill $\blacksquare$

\section{Proof of Theorem \ref{thm: statistical validity}}
\textit{Proof.}
    Define the concentration event
    $$ \mathcal{E}_t = \left\{h(k,x_k^t) - h(k,x_k^*) < \epsilon_t, k\in \calR_t\right\} \text{ and }
    \mathcal{E} = \cap_{t=1}^T\left(\mathcal{E}_t\cap \operatorname{CD}^t\right).$$

    Recall $\calR_t$ is the remaining set of systems at the \textit{beginning} at stage $t$ for $t=1,\ldots,T$. Furthermore, we let   $\calR_{T+1}$ denote the remaining set at the \textit{end} of stage $T$. First notice that, Lemma \ref{lem: prunning validity} indicates $\mathbb{P}(\operatorname{CD}^t|\calR_t,x^t) \ge 1 - \frac{\alpha}{2T}$, which further indicates $\mathbb{P}(\operatorname{CD}^t) \ge 1 - \frac{\alpha}{2T}$ by Tower property. Consequently, with a union bound $\mathcal{E}$ holds with probability at least $1-\alpha$ by Lemma \ref{lem: optimization validity} and Lemma \ref{lem: prunning validity}.
    
    We first show that, under event $\mathcal{E}$, $k^* \in \calR_t$ holds for $t=1,\ldots,T+1$. 
    We prove by induction. 
    First, $k^* \in \calR_1 = \mathcal{K}$. Suppose that $k^*\in\calR_t$ is true for $1\le t\le N$. If $|\calR_t|=1$, then we know $k^*\in\calR_{t'}$ for $t'\ge t$. If $|\calR_t|\ge2$, under $\mathcal{E}_t$, we have 
    $h(k,x_{k}^t) \le h(k,x_{k}^*) + \epsilon_t, \forall k\in\calR_t$. This implies that, for any $k<k^*$, $h(k,x_k^t) - h(k^*,x_{k^*}^t) > h(k,x_k^*) - h(k^*,x_{k^*}^*) -\epsilon_t\ge -\epsilon_t $. 
    Hence, under $\operatorname{CD}^t$, 
    $h(k,x_k^t) - h(k^*,x_{k^*}^t) \le -q_t$ will be declared infeasible, which means $k^*$ will not be removed from the remaining set by comparison with system $k$. Similarly, for $k > k^*$, we can obtain
    $ h(k^*,x_{k^*}^t)- h(k,x_k^t) < \epsilon_t $, which means $h(k^*,x_{k^*}^t) - h(k,x_k^t) \le -q_t$ will be declared feasible. This further implies $k^*$ will not be removed from the remaining set by comparison with $k$. Hence, $k^* \in \calR_{t+1}$. By induction we obtain $k^*\in\calR_t, 1\le t\le T+1$.
    
    Finally, we prove any $k\in\calR_{T+1}$ is an $\epsilon$-optimal system. To see this, if $|\calR_{T+1}| = 1$, then we know $\calR_{T+1} = \{k^*\}$ and $k^*$ is an $\epsilon$-optimal system. Otherwise for any $k \in \calR_{T+1}$,
    $$h(k,x_k^*) - h(k^*,x_{k^*}^*) \le h(k,x_k^T) - h(k^*,x_{k^*}^T) + \epsilon_T \le \epsilon'_T + \epsilon_T = \epsilon, $$
    where the first inequality holds under $\mathcal{E}_N$ by Lemma \ref{lem: optimization validity} and the second inequality holds under $\operatorname{CD}^t$ by Lemma \ref{lem: prunning validity}. This implies any $k\in\calR_{T+1}$ is an $\epsilon$-optimal systemn. This completes the proof.
\hfill $\blacksquare$

\end{document}